\documentclass[a4paper,12pt,twoside]{article}
\usepackage{pst-fill,pst-grad} 
\usepackage{textcomp}
\usepackage[english]{babel}
\usepackage[utf8]{inputenc} 
\usepackage{graphicx}
\usepackage{amsmath}
\usepackage{float}
\usepackage{fancyhdr}
\usepackage[matrix,arrow,curve]{xy}
\usepackage{pstricks} 
\usepackage{amsmath,amsfonts,verbatim,afterpage,theorem,euscript,mathrsfs,amssymb}
\usepackage{amsfonts}
\usepackage{amssymb}
\usepackage{array}
\usepackage{dsfont}
\usepackage{hyperref}

\newtheorem{Definition}{Definition}[section] 
\newtheorem{Proposition}{Proposition}[section]
\newtheorem{Lemme}{Lemma}[section]
\newtheorem{Theoreme}{Theorem}

\def \vu{\vec{u}}

\def \U{\vec{U}}
\def \W{\vec{W}}
\def \V{\vec{V}}
\def \Rt{\mathbb{R}^{3}}

\def \finpv{\hfill $\blacksquare$  \\ \newline }
\def \pv{{\bf{Proof.}}~}

\def \ds{\displaystyle}

\title{\bf Some Liouville theorems for stationary Navier-Stokes equations in Lebesgue and Morrey spaces}

\author{Diego Chamorro\footnote{Laboratoire de Math\'ematiques et Mod\'elisation d'Evry (LaMME) - UMR 8071. Universit\'e d'Evry Val d'Essonne, 23 Boulevard de France, 91037 Evry Cedex, France}$^{\; \;, }$\footnote{email: \text{diego.chamorro@univ-evry.fr}}$\,$, Oscar Jarr\'in$^{*}$ $\,$\footnote{email: \text{oscar.jarrin@dauphine.fr}}$\,$ , Pierre-Gilles Lemari\'e-Rieusset$^{*,}$\footnote{email: \text{plemarie@univ-evry.fr}}} 
\begin{document}
\maketitle

\begin{scriptsize}
\abstract{Uniqueness of Leray solutions of the 3D Navier-Stokes equations is a challenging open problem. In this article we will study this problem for the 3D stationary Navier-Stokes equations and under some additional hypotheses, stated in terms of Lebesgue and Morrey spaces, we will show that the trivial solution $\U=0$  is the unique solution. This type of results are known as Liouville theorems.}\\[3mm]
\end{scriptsize}
\section{Introduction} 
In this article we study uniqueness of weak solutions to the stationary and incompressible Navier-Stokes equations in the whole space $\Rt$:
\begin{equation}\label{N-S-f-nulle-intro}
-\Delta \U+(\U \cdot \vec{\nabla})\U +\vec{\nabla}P = 0, \qquad div(\U)=0, \\
\end{equation}
where $\U:\Rt\longrightarrow \Rt$ is the velocity and $P:\Rt\longrightarrow \mathbb{R}$ is the pressure. Recall that a weak solution of equations (\ref{N-S-f-nulle-intro}) is a couple $(\U,P)\in L^{2}_{loc}(\Rt)\times \mathcal{D}^{'}(\Rt)$ which verifies these equations in the distributional sense. Recall also that we can concentrate our study in the velocity $\U$ since we have the identity $P=\frac{1}{(-\Delta)}div\big((\U \cdot \vec{\nabla})\U\big )$.\\

It is clear that the trivial solution $\U=0$ satisfies (\ref{N-S-f-nulle-intro}) and it is natural to ask if this is the unique solution of these equations. In the general setting of the space $L^{2}_{loc}(\Rt)$, the answer is negative: indeed, if we define the function $\psi:\Rt \longrightarrow \mathbb{R}$ by $\psi(x_1,x_2,x_3)=\frac{x^{2}_{1}}{2}+\frac{x^{2}_{2}}{2}-x^{2}_{3}$ and if we set the functions $\U$ and $P$ by  the identities
$$\U(x_1,x_2,x_3)=\vec{\nabla}\psi(x_1,x_2,x_3)=(x_1, x_2, -2x_3),\quad \mbox{and}\quad P(x_1,x_2,x_3)=-\frac{1}{2}\vert \U(x_1,x_2,x_3)\vert^2,$$
then we have $\U \in L^{2}_{loc}(\Rt)$ (since $\vert \U(x)\vert \approx \vert x \vert$) and using basic rules of vector calculus we have that the couple $(\U,P)$ given by the expressions above satisfies (\ref{N-S-f-nulle-intro}).\\

Thus, due to this lack of uniqueness in the general setting of space $L^{2}_{loc}(\Rt)$ we are interested in the following problem (also known as \emph{Liouville problem}): find a functional space $E \subset L^{2}_{loc}(\Rt)$ such that if $\U \in L^{2}_{loc}(\Rt)$ is a solution of equations (\ref{N-S-f-nulle-intro}) and if $\U \in E$, then $\U=0$. \\

A well-known result on the Liouville problem for equation (\ref{N-S-f-nulle-intro})  is given  in the  book  \cite{Galdi} of G. Galdi where it is shown that to prove the identity $\U=0$, we need a certain decrease at infinity of the solution. More precisely, if the solution $\U \in L^{2}_{loc}(\Rt)$ verifies the additional hypothesis $\U \in L^{\frac{9}{2}}(\Rt)$ then we have $\U=0$ (see \cite{Galdi}, Theorem X.9.5, page $729$). This result has been improved in different settings: D. Chae and J. Wolf gave a logarithmic improvement of Galdi’s result in \cite{ChaeWolf}. Moreover, H. Kozono \emph{et.al.} prove in \cite{Kozono} that $\U=0$ when $\U \in L^{\frac{9}{2},\infty}(\Rt)$ and with additional conditions on the decay (in space variable)  of the vorticity $\vec{w}= \vec{\nabla}\wedge \U$. For more references on the Liouville problem for the stationary Navier-Stokes equations see also the articles  \cite{ChaeYoneda}, \cite{ChaeWeng} and \cite{Koch} and the references therein. \\

Another interesting result was given by G. Seregin in \cite{Ser2015} where the hypothesis $\U \in L^{\frac{9}{2}}(\Rt)$ is replaced by the condition $\U \in L^{6}(\Rt)\cap BMO^{-1}(\Rt)$: here the solution $\U$ decrease slowly to infinity since we only have $ \U \in L^{6}(\Rt)$ and thus the extra hypothesis $BMO^{-1}(\Rt)$ is added to get the desired identity $\U=0$. \\

In our first theorem we generalize previous results and we study the Liouville problem in the setting of Lebesgue spaces:
\begin{Theoreme}\label{Theo:Lebesgue} Let $\U \in L^{2}_{loc}(\Rt)$ be a weak solution of the stationary Navier-Stokes  equations (\ref{N-S-f-nulle-intro}). 
\begin{enumerate}
\item[1)] If $\U \in L^p(\Rt)$ with $3 \leq p \leq \frac{9}{2}$, then $\U=0$.
\item[2)] If $\U \in L^p(\Rt)\cap \dot{B}^{\frac{3}{p}-\frac{3}{2}, \infty}_{\infty}(\Rt)$ with $\frac{9}{2}<p<6$, then $\U=0$.
\end{enumerate}
\end{Theoreme}
In the second point above, since $\frac{3}{p}-\frac{3}{2}<0$ we can characterize the Besov space $\dot{B}^{\frac{3}{p}-\frac{3}{2},\infty}_{\infty}(\Rt)$ as the set of distributions $f\in \mathcal{S}'(\Rt)$ such that $\ds{ \Vert f \Vert_{\dot{B}^{\frac{3}{p}-\frac{3}{2},\infty}_{\infty}}=\sup_{t>0} t^{\frac{1}{2}(\frac{3}{2}-\frac{3}{p})} \Vert h_{t}\ast f \Vert_{L^{\infty}}}<+\infty$ where $h_t$ denotes the heat kernel.\\

It is worth noting here that the space $L^{\frac{9}{2}}(\Rt)$ seems to be a \emph{limit space} to solve the Liouville problem in the sense that if $3 \leq p \leq \frac{9}{2}$ we do not need any extra information, but if  $\frac{9}{2}<p<6$ we need an additional hypothesis given in terms of Besov spaces. Remark also that, to the best of our knowledge, the Liouville problem for stationary Navier-Stokes equations in the Lebesgue spaces $L^{p}(\Rt)$ with $1\leq p <3$ or $6\leq p\leq +\infty$ is still an open problem.\\

More recently G. Seregin \cite{Ser2016} replaced the hypothesis $\U \in L^6(\Rt) \cap BMO^{-1}(\Rt)$ by a couple of homogeneous Morrey spaces $\dot{M}^{p,q}(\Rt)$. Recall that for $1<p\leq q<+\infty$ the space $\dot{M}^{p,q}(\Rt)$ is defined as the functions $f\in L^{p}_{loc}(\Rt)$ such that
\begin{equation}\label{def-morrey} 
\ds{\Vert f \Vert_{\dot{M}^{p,q}}= \sup_{x_0\in \Rt, \,r>0} \left(  r^{\frac{3}{q}-\frac{3}{p}} \times   \left(\int_{B(x_0,r)} \vert f(x)\vert^p dx \right)^{\frac{1}{p}}\right) <+\infty}. 
\end{equation}
This space is an homogeneous space of degree $-\frac{3}{q}$ and in Theorem 1.1 of \cite{Ser2016}  it is shown that if the solution $\U \in L^{2}_{loc}(\Rt)$ verifies $\U \in  \dot{M}^{2,6}(\Rt) \cap \dot{M}^{\frac{3}{2},3}(\Rt)$ then we have $\U=0$. \\

If we compare the condition $\U \in L^6(\Rt)$ and $BMO^{-1}(\Rt)$ given in \cite{Ser2015} with the hypothesis $\U \in \dot{M}^{2,6}(\Rt)\cap \dot{M}^{\frac{3}{2},3}(\Rt)$ given in \cite{Ser2016}, we can observe that the shift to Morrey spaces preserves the homogeneity: $L^6(\Rt)$ is substituted by the Morrey space $\dot{M}^{2,6}(\Rt)$ with the same homogeneous degree $-1$ while $BMO^{-1}(\Rt)$ is replaced by the Morrey space $\dot{M}^{\frac{3}{2},3}(\Rt)$, also with homogeneous degree $-1$. \\

Following these ideas we study the Liouville problem in the setting of Morrey spaces for equations (\ref{N-S-f-nulle-intro}) and we generalize the result obtained in \cite{Ser2016} in the following way: 

\begin{Theoreme}\label{Theo:Seregin-generalise} Let $\U \in L^{2}_{loc}(\Rt)$ be a weak solution of the stationary Navier-Stokes equations (\ref{N-S-f-nulle-intro}). If $\U\in \dot{M}^{2,3}(\Rt)\cap \dot{M}^{2,q}(\Rt)$ with $3<q<+\infty$, then we have $\U=0$. 
\end{Theoreme}

We observe here that we kept an homogeneous  Morrey space of degree $-1$, namely $\dot{M}^{2,3}(\Rt)$, but the space $\dot{M}^{2,6}(\Rt)$ used previously in \cite{Ser2016} is now replaced by \emph{any} Morrey space $\dot{M}^{2,q}(\Rt)$  which is an homogeneous space of degree $-1<-\frac{3}{q}<0$.\\

A natural question raises: it is possible to consider a single Morrey space in order to solve the Liouville problem for equation (\ref{N-S-f-nulle-intro})? The answer is positive, but we need to introduce the following functional space. 
\begin{Definition}\label{Def:sous-espace-morrey} Let  $1<p\leq q <+\infty$. We define the space $\overline{M}^{p,q}(\Rt)$ as the closure of the test functions $\mathcal{C}^{\infty}_{0}(\Rt)$ in the Morrey space $\dot{M}^{p,q}(\Rt)$. 
\end{Definition}

The space $\overline{M}^{p,q}(\Rt)$ is of course smaller than $\dot{M}^{p,q}(\Rt)$, and for suitable values of the parameters $p,q$ we have the following result.
\begin{Theoreme}\label{Theo:Morrey-2} Let $2<p\leq 3$ and consider the space $\overline{M}^{p,3}(\Rt)$  given by Definition \ref{Def:sous-espace-morrey} above. Let $\U \in L^{2}_{loc}(\Rt)$ be a weak solution of the stationary Navier-Stokes equations  (\ref{N-S-f-nulle-intro}). If  $\U\in \overline{M}^{p,3}(\Rt)$  then $\U=0$. 
\end{Theoreme}

The reason why  we prove the uniqueness of the solution $\U=0$ in the setting of the space $\overline{M}^{p,3}(\Rt)$ and not in the more general setting of the space $\dot{M}^{p,3}(\Rt)$ is purely technical as we will explain in details in Section \ref{Sec:Theo:Morrey-2}. \\

This article is organized as follows: in Section \ref{sec:Lebesgue} we study the Liouville problem for equations (\ref{N-S-f-nulle-intro}) in the setting of Lebesgue space. Then, in Section \ref{Ths-Morrey} we study the Liouville problem in the setting of Morrey spaces where we prove Theorem \ref{Theo:Seregin-generalise} and Theorem \ref{Theo:Morrey-2}. Section \ref{appendix} is reserved for a  technical lemma. 

\section{The Liouville problem in Lebesgue spaces}\label{sec:Lebesgue}
We prove here Theorem \ref{Theo:Lebesgue} and from now on $\U\in L^{2}_{loc}(\Rt)$ will be a weak solution of the stationary Navier-Stokes equations (\ref{N-S-f-nulle-intro}).
\begin{enumerate}
\item[1)]  Assume that $\U \in L^p(\Rt)$ with  $3\leq p \leq \frac{9}{2}$. We are going to prove the identity $\U=0$ and for this we will follow the main ideas of \cite{Galdi} (Theorem X.9.5, page $729$). We start then by introducing the following cut-off function: let  $\theta \in \mathcal{C}^{\infty}_{0}(\Rt)$ be  such that  $0\leq \theta \leq 1$, $\theta(x)=1$ if $\vert x \vert <\frac{1}{2}$ and  $\theta(x)=0$ if $\vert x \vert \geq 1$. Let now $R>1$ and define the function $\theta_R(x)=\theta\left( \frac{x}{R}\right)$, we have then $\theta_R(x)=1$ if $\vert x \vert < \frac{R}{2}$ and $\theta_R(x)=0$ if $\vert x \vert \geq R$.\\

Now, we multiply equation (\ref{N-S-f-nulle-intro}) by the function $\theta_R \U$, then  we integrate on the ball $B_R=\{x \in \Rt: \vert x \vert <R \}$ to obtain the following identity
$$\int_{B_R} \left(- \Delta \U +(\U \cdot \vec{\nabla})\U +\vec{\nabla}P \right)\cdot (\theta_R \U) dx=0.$$
Observe that, since $\U \in L^p(\Rt)$ with $3\leq p\leq \frac{9}{3}$ then $\U \in L^{3}_{loc}(\Rt)$ and by Theorem X.1.1 of the book \cite{Galdi} (page $658$), we have  $\U \in \mathcal{C}^{\infty}(\Rt)$ and $P\in \mathcal{C}^{\infty}(\Rt)$. Thus, all the terms in the identity above are well-defined  and we have
\begin{equation}\label{iden-1}
\int_{B_R} - \Delta \U\cdot \left( \theta_R \U\right) +(\U \cdot \vec{\nabla})\U \cdot \left(  \theta_R \U\right) +\vec{\nabla}P\cdot \left( \theta_R \U\right) dx=0.
\end{equation} 
We study now each term in this identity.  For the first term in (\ref{iden-1}), integrating by parts and since  $\theta_R(x)=0$ if $\vert x \vert \geq R$, then we write

\begin{eqnarray} \nonumber 
\int_{B_R} - \Delta \U\cdot \left( \theta_R \U\right) dx &= &-\sum_{i,j=1}^{3} \int_{B_R} (\partial^{2}_{j} U_i) ( \theta_R U_i) dx= \sum_{i,j=1}^{3} \int_{B_R} \partial_{j}U_i \partial_{j}(\theta_R U_i) dx \\  \nonumber
&=&  \sum_{i,j=1}^{3} \int_{B_R} (\partial_j U_i) (\partial_j \theta_R)    U_i dx +
 \sum_{i,j=1}^{3} \int_{B_R} (\partial_j U_i)  \theta_R (  \partial_j U_i) dx\\ \nonumber
 &=&  \sum_{i,j=1}^{3}   \int_{B_R} (\partial_j \theta_R) (\partial_j U_i)    U_i dx +  \sum_{i,j=1}^{3} \int_{B_R} \theta_R  (\partial_j U_i)^2  dx\\ \nonumber
& =&  \sum_{i,j=1}^{3} \int_{B_R} (\partial_j \theta_R )  \partial_{j} \left( \frac{U^{2}_{i}}{2}\right)dx + \int_{B_R}\theta_R \vert \vec{\nabla}\otimes U \vert^{2} dx\\ 
&=& - \int_{B_R} \Delta \theta_R  \left( \frac{\vert U\vert^2}{2}  \right)dx +   \int_{B_R}\theta_R \vert \vec{\nabla}\otimes U \vert^{2} dx.  \label{ipp1}
\end{eqnarray}

For the second term in (\ref{iden-1}) we write

\begin{eqnarray} \nonumber 
\int_{B_R} (\U \cdot \vec{\nabla})\U\cdot ( \theta_R \U )dx &=& \sum_{i,j=1}^{3} \int_{B_R} U_j (\partial_j U_i) ( \theta_R U_i ) dx =  \sum_{i,j=1}^{3} \int_{B_R} \theta_R U_j  (\partial_j U_i) U_i dx\\
&=&  \sum_{i,j=1}^{3} \int_{B_R} \theta_R U_j  (\partial_j \left(\frac{U^{2}_{i}}{2}\right)) dx,   \label{ipp2}
\end{eqnarray}
but, as $div(\U)=0$ and then integrating by parts  we can write 
\begin{equation} \label{ipp3}
\sum_{i,j=1}^{3} \int_{B_R} \theta_R U_j  (\partial_j \left(\frac{U^{2}_{i}}{2}\right)) dx =\sum_{i,j=1}^{3} \int_{B_R} \theta_R \partial_{j} \left( U_j \frac{U^{2}_{i}}{2} \right)dx  - \int_{B_R} \vec{\nabla} \theta_R \cdot \left( \frac{\vert \U\vert^2}{2} \U\right)dx.  
\end{equation}
For the third term in (\ref{iden-1}), integrating by parts and  since $div(\U)=0$ then we have 
\begin{eqnarray} \nonumber
\int_{B_R} \vec{\nabla}P\cdot ( \theta_R \U )dx&=& 	\sum_{i=1}^{3}\int_{B_R}(\partial_i P) \theta_R U_i dx=- \sum_{i=1}^{3}\int_{B_R} P \partial_i (\theta_R U_i)dx \\
&=& - \sum_{i=1}^{3}\int_{B_R} P (\partial_i \theta_R)  ( U_i) dx = - \int_{B_R}\vec{\nabla}\theta_R \cdot(P \U) dx. \label{ipp4}
	\end{eqnarray}

With these identities  and getting back to  equation (\ref{iden-1}) we can write 
\begin{equation*}
- \int_{B_R} \Delta \theta_R  \left( \frac{\vert U\vert^2}{2}  \right)dx +   \int_{B_R}\theta_R \vert \vec{\nabla}\otimes U \vert^{2} dx - 	\int_{B_R} \vec{\nabla} \theta_R \cdot \left( \frac{\vert \U\vert^2}{2} \U\right)dx - \int_{B_R}\vec{\nabla}\theta_R \cdot(P \U) dx=0,
\end{equation*}
hence we get
\begin{equation}\label{Identity1}
\int_{B_R} \theta_R \vert \vec{\nabla} \otimes \U \vert^2 dx  = \int_{B_R} \Delta \theta_R  \frac{\vert \U \vert^2}{2} dx + \int_{B_R} \vec{\nabla} \theta_R \cdot \left( \left( \frac{\vert \U \vert^2}{2}+P \right) \U \right) dx.
\end{equation}
On the other hand,  as $\theta_R (x)=1$ if $\vert x \vert < \frac{R}{2}$ then we have 
$$\ds{ \int_{B_{\frac{R}{2}}} \vert \vec{\nabla}\otimes \U \vert^2 dx \leq  \int_{B_R} \theta_R \vert \vec{\nabla} \otimes \U \vert^2 dx},$$ 
and by identity (\ref{Identity1}) we obtain
\begin{eqnarray}\label{estim-aux-Lioville-1}\nonumber
 \int_{B_{\frac{R}{2}}} \vert \vec{\nabla}\otimes \U \vert^2 dx & \leq &\int_{B_R} \Delta \theta_R  \frac{\vert \U \vert^2}{2} dx + \int_{B_R} \vec{\nabla} \theta_R \cdot \left( \left( \frac{\vert \U \vert^2}{2}+P \right) \U \right) dx \\	
 &\leq & I_1(R) +I_2(R), 
\end{eqnarray}  
and we will prove  that  $\ds{\lim_{R \longrightarrow +\infty}I_i(R)=0}$ for $i=1,2$. \\

Indeed, for the term $I_1(R)$, by H\"older inequalities (with $\frac{1}{q}+\frac{2}{p}=1$) we have
$$ I_1(R)\leq \left( \int_{B_R} \vert \Delta \theta_R \vert^q dx\right)^{\frac{1}{q}} \left( \int_{B_R} \vert \U \vert^p dx\right)^{\frac{2}{p}} \leq \left( \int_{B_R} \vert \Delta \theta_R \vert^q dx\right)^{\frac{1}{q}} \Vert \U \Vert^{2}_{L^p}.$$ Moreover, as $\theta_R(x)=\theta\left(\frac{x}{R}\right)$ we have 
$ \ds{ \left( \int_{B_R} \vert \Delta \theta_R \vert^q dx\right)^{\frac{1}{q}}= R^{\frac{3}{q}-2} \times \Vert  \Delta \theta \Vert_{L^q(B_1)}}$,  
and as $\frac{1}{q}+\frac{2}{p}=1$ then we can write 
$ \ds{ I_1(R)\leq R^{1-\frac{6}{p}} \times  \Vert \Delta \theta \Vert_{L^q(B_1)} \Vert \U \Vert^{2}_{L^p}}$. \\
\\ 
In this estimate we observe that  since $3\leq p\leq \frac{9}{2}$ then $-1\leq 1-\frac{6}{p}\leq-\frac{1}{3}$ and thus we get $\ds{\lim_{R\longrightarrow +\infty} I_1(R)=0}$. \\

We study now the term $I_2(R)$ in (\ref{estim-aux-Lioville-1}). Recall that $\theta_R (x)=1$ if $\vert x \vert <\frac{R}{2}$ and $\theta_R(x)=0$ if $\vert x \vert \geq R$, so we have  $supp \left( \vec{\nabla}\theta_R \right) \subset \{ x\in \Rt: \frac{R}{2} < \vert x \vert <R \}=\mathcal{C}(\frac R2, R)$ and we can write 
\begin{equation*}
I_2(R) = \int_{B_R} \vec{\nabla} \theta_R \cdot \left( \left( \frac{\vert \U \vert^2}{2}+P \right) \U \right) dx=\int_{\mathcal{C}(\frac R2, R)} \vec{\nabla} \theta_R \cdot \left( \left( \frac{\vert \U \vert^2}{2}+P \right) \U \right) dx,
\end{equation*}
hence we have
\begin{eqnarray*} 
|I_2(R)| &\leq & \frac{1}{2}\int_{\mathcal{C}(\frac R2, R)} \vert \vec{\nabla} \theta_R \vert \vert \U\vert^3 dx + \int_{\mathcal{C}(\frac R2, R)} \vert \vec{\nabla} \theta_R \vert \vert P \vert \vert \U\vert dx\\
&\leq & (I_2)_a(R)+(I_2)_b(R),	
\end{eqnarray*}	
and we will prove now that  $\ds{ \lim_{R\longrightarrow+\infty} (I_2)_a(R)=0}$ and $ \ds{\lim_{R \longrightarrow +\infty} (I_2)_b(R)=0}$. \\
\\
For the term $(I_2)_a(R)$,  by H\"older inequalities (with $\frac{1}{r}+\frac{3}{p}=1$)  we have
\begin{equation}\label{eq-aux1} 
(I_2)_a(R) \leq  \left( \int_{\mathcal{C}(\frac R2, R)}\vert \vec{\nabla} \theta_R \vert^{r} dx \right)^{\frac{1}{r}} \left( \int_{\mathcal{C}(\frac R2, R)} \vert \U \vert^p dx\right)^{\frac{3}{p}},
\end{equation}   and we study now the first term in the right side. As $\theta_R(x)=\theta\left(\frac{x}{R}\right)$  then we have $\ds{\left( \int_{\mathcal{C}(\frac R2, R)}\vert \vec{\nabla} \theta_R \vert^{r} dx \right)^{\frac{1}{r}} \leq R^{\frac{3}{r}-1}}\Vert \vec{\nabla}\theta \Vert_{L^r}$, and since $\frac{1}{r}=1-\frac{3}{p}$ then we have $\frac{3}{r}-1=2-\frac{9}{p}$ and thus we write $\ds{\left( \int_{\mathcal{C}(\frac R2, R)}\vert \vec{\nabla} \theta_R \vert^{r} dx \right)^{\frac{1}{r}} \leq R^{2-\frac{9}{p}}\Vert \theta \Vert_{L^r}}$. But, since $3\leq p\leq \frac{9}{2}$ then we have $-1\leq 2-\frac{9}{p}\leq 0$, and  since  $R>1$ then we get  $R^{2-\frac{9}{p}}\leq 1$. So,  by the last inequality we can write  
 \begin{equation}\label{estim-aux1}
 \left( \int_{\mathcal{C}(\frac R2, R)}\vert \vec{\nabla} \theta_R \vert^{r} dx \right)^{\frac{1}{r}} \leq \Vert \vec{\nabla} \theta \Vert_{L^r}. 
 \end{equation}
 With this estimate and getting back to estimate (\ref{eq-aux1}) we can write  $$\ds{(I_2)_a(R) \leq  \Vert \vec{\nabla} \theta \Vert_{L^r}  \Vert \U \Vert^{3}_{L^p(\mathcal{C}(\frac R2, R))}},$$ 
and since $\U \in L^p(\Rt)$  then we have $\ds{\lim_{R \longrightarrow +\infty} \Vert \U \Vert_{L^p(\mathcal{C}\left(\frac{R}{2},R\right))}=0}$ and we obtain $$\ds{\lim_{R\longrightarrow+\infty} (I_2)_a(R)=0.}$$ 
\\
For the term  $ (I_2)_b(R)$,  by   H\"older inequalities  (with $\frac{1}{r}+\frac{3}{p}=1$) and  by estimate (\ref{estim-aux1})  we can write
\begin{eqnarray}\label{estim-aux2} \nonumber
	(I_2)_b (R)&\leq & \int_{\mathcal{C}(\frac R2, R)} \vert \vec{\nabla} \theta_R \vert \vert P \vert \vert \U\vert dx \leq \left( \int_{\mathcal{C}(\frac R2, R)} \vert \vec{\nabla} \theta_R \vert^{r}dx \right)^{\frac{1}{r}} \left( \int_{\mathcal{C}(\frac R2, R)} ( \vert P \vert \vert \U \vert)^{\frac{p}{3}} dx \right)^{\frac{3}{p}}\\
	&\leq & \Vert \vec{\nabla} \theta \Vert_{L^r} \left( \int_{\mathcal{C}(\frac R2, R)} ( \vert P \vert \vert \U \vert)^{\frac{p}{3}} dx \right)^{\frac{3}{p}}.
\end{eqnarray}
But, recall  that since the velocity $\U$ belongs to the space $L^p(\Rt)$ then pressure $P$ belongs to the space $L^{\frac{p}{2}}(\Rt)$. Indeed, we write
\begin{equation}\label{pression}
\ds{P=\sum_{i,j=1}^{3}\frac{1}{-\Delta}\partial_i \partial_j (U_i U_j)=\sum_{i,j=1}^{3}\mathcal{R}_i\mathcal{R}_j (U_i U_j)},
\end{equation}  where $\mathcal{R}_i=\frac{\partial_i}{\sqrt{-\Delta}}$
denotes the i-th Riesz transform. By the continuity of the operator  $\mathcal{R}_i\mathcal{R}_j$ on Lebesgue  spaces $L^q(\Rt)$ (with $1<q<+\infty$)  and applying the H\"older inequalities we get $\Vert P \Vert_{L^{\frac{p}{2}}}\leq c \Vert \U \Vert^{2}_{L^p}$. \\
\\
Then, getting back to estimate (\ref{estim-aux2}), always by   H\"older inequalities  (with $\frac{2}{p}+\frac{1}{p}=\frac{3}{p}$) we write 
$$ (I_2)_b(R) \leq   \Vert \vec{\nabla} \theta \Vert_{L^r} \left( \int_{\mathcal{C}(\frac R2, R)} \vert P \vert^{\frac{p}{2}} dx \right)^{\frac{2}{p}} \left( \int_{\mathcal{C}(\frac R2, R)} \vert \U \vert^{p}dx \right)^{\frac{1}{p}},$$
and since $\U \in L^{p}(\Rt)$ and $P\in L^{\frac{p}{2}}(\Rt)$ then we get  $\ds{\lim_{R\longrightarrow+\infty} (I_2)_b(R)=0}$. We have proven that $\ds{\lim_{R\longrightarrow+\infty} I_2(R)=0}$. \\
\\
Now with the information $\ds{\lim_{R\longrightarrow+\infty} I_i(R)=0}$ for $i=1,2$ we get back to estimate (\ref{estim-aux-Lioville-1}) and we can deduce that $ \ds{\int_{\Rt} \vert \vec{\nabla} \otimes \U \vert^2 dx=0}$. But, recall that  by the Hardy-Littlewood-Sobolev inequalities  we have $\Vert \U \Vert_{L^{6}(\Rt)}\leq c \Vert \U \Vert_{\dot{H}^1(\Rt)}$ and thus  we have the identity  $\U=0$. \\ 
\item[2)] We suppose now  $\U \in L^p(\Rt)\cap \dot{B}^{\frac{3}{p}-\frac{3}{2},\infty}_{\infty}(\Rt)$ with  $\frac{9}{2}<p<6$ and we will prove that  $\U=0$. For this we will follow some ideas of the article \cite{Ser2015} and the first thing to do is to prove the following proposition.
\begin{Proposition}\label{Prop:Lp-Besov} Let $\frac{9}{2}<p<6$ and let $\U \in L^p\cap \dot{B}^{\frac{3}{p}-\frac{3}{2},\infty}_{\infty}(\Rt)$ be a weak solution of the stationary Navier-Stokes equations  (\ref{N-S-f-nulle-intro}). Then  $\U\in \dot{H}^1(\Rt)$ and we have $ \Vert \U \Vert_{\dot{H}^1}\leq c\left(1+ \Vert \U \Vert_{\dot{B}^{\frac{3}{p}-\frac{3}{2},\infty}_{\infty}}\right)  \Vert \U \Vert_{L^p}$.  
\end{Proposition}
\textbf{Proof.} To prove this result we need to verify the following estimate (also called a \emph{Cacciopoli type inequality} \cite{Ser2015}, \cite{Ser2016}): let $R>1$ and let the ball $B_R=\lbrace x\in \Rt: \vert x \vert <R \rbrace$, then we have
\begin{equation}\label{Cacciopoli}
\int_{B_{\frac{R}{2}}}\vert \vec{\nabla}\otimes \U (x)\vert^2dx \leq C(\U, R)  \Vert \U\Vert^{2}_{L^p},
\end{equation} 
where  $\ds{C(\U,R)=c \left(R^{1-\frac{6}{p}}+1 \right)  \times  \left(1+ \Vert \U \Vert_{\dot{B}^{\frac{3}{p}-\frac{3}{2},\infty}_{\infty}}\right)^{2}}$, 
and where $c>0$ is a  constant which does not depend of the  solution $\U$ nor  of $R>1$. \\
\\
To verify  (\ref{Cacciopoli}) we start by introducing the test functions $\varphi_R$ and $\W_R$ as follows:  for a fixed $R>1$, we define first the function $\varphi_R\in \mathcal{C}^{\infty}_{0}(\Rt)$ by $0\leq \varphi_R\leq 1$ such that  for $\frac{R}{2}\leq \rho<r<R$ we have  $\varphi_R(x)=1$ if $\vert x \vert<\rho$,  $ \varphi_R(x)=0$  if $\vert x \vert\geq r$ and 
\begin{equation}\label{control-test}
\Vert \vec{\nabla}\varphi_R\Vert_{L^{\infty}}\leq \frac{c}{r-\rho}.
\end{equation}Next we define the function $\W_R$ as the solution of the problem
\begin{equation}\label{eq_W_R}
div(\W_R)=\vec{\nabla}\varphi_R\cdot \U, \quad \text{over}\,\, B_r, \quad \text{and}\quad \W_R=0 \,\, \text{over}\,\,  \partial B_r,
\end{equation} 
where $\partial B_r=\{x\in \Rt: \vert x \vert=r \}$.  Existence of such function $\W_R$ is assured by Lemma  $III. 3.1$ (page 162) of the book  \cite{Galdi} and where it is proven that $\W_R\in W^{1,p}(B_r)$ with
\begin{equation}\label{estim-W}
\Vert \vec{\nabla}\otimes \W_R\Vert_{L^p(B_r)}\leq c \Vert   \vec{\nabla}\varphi_R\cdot \U \Vert_{L^p(B_r)}.
\end{equation} 
Once we have defined the functions $\varphi_R$ and $\W_R$ above, we consider now the function $\varphi_R \U-\W_R$ and we write  
\begin{equation}\label{iden1}
\int_{B_r} \left( -\Delta \U +(\U \cdot \vec{\nabla})\U +\vec{\nabla}P\right)\cdot \left( \varphi_R \U-\W_R \right)dx=0.
\end{equation}
Remark that since $\U \in L^p(\Rt)$ with $\frac{9}{2}<p<6$ then $\U \in L^{3}_{loc}(\Rt)$ and always by Theorem X.1.1  of the book \cite{Galdi} (page $658$) we have $\U \in \mathcal{C}^{\infty}(\Rt)$ and $P\in \mathcal{C}^{\infty}(\Rt)$ and thus every term in the last identity  is well-defined.\\

In the identity (\ref{iden1}), we start by studying the third term   $\displaystyle{\int_{B_r}\vec{\nabla}P\cdot \left( \varphi_R \U-\W_R \right)dx}$ and by an integration by parts we write
$$ \int_{B_r} \vec{\nabla}P \cdot \left( \varphi_R \U-\W_R \right)dx=-\int_{B_r}P\left( \vec{\nabla}\varphi_R\cdot \U+ \varphi_R\, div(\U)-div(\W_R)\right)dx,$$
but since $\W_R$ is a solution of problem (\ref{eq_W_R}) and since $div(\U)=0$ then we can write  $\displaystyle{\int_{B_r}\vec{\nabla}P\cdot \left( \varphi_R \U-\W_R \right)dx=0}$ and thus identity (\ref{iden1}) can be written as:
\begin{equation}\label{iden2}
\int_{B_r} -\Delta \U \cdot \left( \varphi_R \U-\W_R \right)dx+ \int_{B_r} \left( (\U \cdot \vec{\nabla})\U \right)\cdot \left( \varphi_R \U-\W_R \right)dx=0.
\end{equation}
In this equation above we study now the   term  $\ds{\int_{B_r} -\Delta \U \cdot \left( \varphi_R \U-\W_R \right)dx}$ and  always integrating by parts we have
\begin{eqnarray*}
 & & \int_{B_r} -\Delta \U \cdot \left( \varphi_R \U-\W_R \right)dx = \sum_{i,j=1}^{3} \int_{B_r} ( \partial_{j} U_i) \partial_{j} \left( \varphi_R U_i - (W_{R})_i \right)dx  \\ \nonumber 
&=& \sum_{i,j=1}^{3} \int_{B_r} \partial_j U_i (\partial_j \varphi_R)U_i dx + \sum_{i,j=1}^{3} \int_{B_r} \varphi_R ( \partial_j U_i )^2 dx -  \sum_{i,j=1}^{3} \int_{B_r} (\partial_j U_i) \partial_j (W_R)_i dx\\
&=& \sum_{i,j=1}^{3} \int_{B_r} \partial_j U_i (\partial_j \varphi_R)U_i dx + \int_{B_r} \varphi_R \vert \vec{\nabla}\otimes \U \vert^2-\sum_{i,j=1}^{3} \int_{B_r} (\partial_j U_i) \partial_j (W_R)_i dx.
\end{eqnarray*}  
With this identity we get back to equation (\ref{iden2}) and we can write 
\begin{eqnarray*}
& &\sum_{i,j=1}^{3} \int_{B_r} \partial_j U_i (\partial_j \varphi_R)U_i dx + \int_{B_r} \varphi_R \vert \vec{\nabla}\otimes \U \vert^2-\sum_{i,j=1}^{3} \int_{B_r} (\partial_j U_i) \partial_j (W_R)_i dx 	\\
&+& \int_{B_r} \left( (\U \cdot \vec{\nabla})\U \right)\cdot \left( \varphi_R \U-\W_R \right)dx=0,
\end{eqnarray*}
hence we have
\begin{eqnarray} \nonumber
\int_{B_r}\varphi_R \vert \vec{\nabla}\otimes \U\vert^2dx &=& - \sum_{i,j=1}^{3} \int_{B_r} \partial_j U_i (\partial_j \varphi_R)U_i dx +  \sum_{i,j=1}^{3} \int_{B_r} (\partial_j U_i) \partial_j (W_R)_i dx \\
& &-\int_{B_r}\left( (\U\cdot \vec{\nabla})\cdot \U \right)\cdot (\varphi_R \U-\W_R)dx\nonumber\\
&=&I_1+I_2+I_3. \label{identity_ipp} 
\end{eqnarray}
Now, we must study the terms  $I_1,I_2$ and $I_3$ above and for this  we decompose our study in two  technical lemmas:   

\begin{Lemme}\label{Lemme1} Let $\frac{9}{2}<p<6$ and let $\U \in L^p\cap \dot{B}^{\frac{3}{p}-\frac{3}{2},\infty}_{\infty}(\Rt)$ be a weak solution of the stationary Navier-Stokes equations  (\ref{N-S-f-nulle-intro}). Then there exists a  constant  $c>0$ (which does not depend of $R,r,\rho$ and  $\U$ ) such that  
$$|I_1|+|I_2|\leq   c\frac{R^{3\left( \frac{1}{2}-\frac{1}{p}\right)}}{r-\rho}\left(\int_{B_r}\vert \vec{\nabla}\otimes \U\vert^2dx\right)^{\frac{1}{2}} \left(\int_{B_r}\vert  \U\vert^pdx\right)^{\frac{1}{p}}.$$
\end{Lemme}
 \textbf{Proof.} For the term $I_1$ in identity (\ref{identity_ipp}),  by  the Cauchy-Schwarz inequality we write 
$$ \vert I_1 \vert= \left| \sum_{i,j=1}^{3} \int_{B_r} \partial_j U_i (\partial_j \varphi_R)U_i dx \right| \leq  c \left(\int_{B_r}\vert \vec{\nabla}\otimes \U\vert^2dx\right)^{\frac{1}{2}}\left(\int_{B_r}\vert \vec{\nabla}\varphi_R\otimes \U\vert^2dx\right)^{\frac{1}{2}} ,$$ then in the second term of the quantity in the right side   we apply the  H\"older inequalities  (with $\frac{1}{2}=\frac{1}{q}+\frac{1}{p}$) and since $\Vert \vec{\nabla}\varphi_R\Vert_{L^{\infty}}\leq \frac{c}{r-\rho}$ we can write
$$ \left(\int_{B_r}\vert \vec{\nabla}\varphi_R\otimes \U\vert^2dx\right)^{\frac{1}{2}} \leq  \left(\int_{B_r}\vert \vec{\nabla}\varphi_R\vert^qdx\right)^{\frac{1}{q}}   \left(\int_{B_R}\vert  \U\vert^pdx\right)^{\frac{1}{p}} \leq c\frac{r^{\frac{3}{q}}}{r-\rho}\left(\int_{B_r}\vert  \U\vert^pdx\right)^{\frac{1}{p}},$$ and thus we have the estimate
$$  \vert I_1 \vert \leq c\frac{r^{\frac{3}{q}}}{r-\rho} \left(\int_{B_r}\vert \vec{\nabla}\otimes \U\vert^2dx\right)^{\frac{1}{2}} \left(\int_{B_r}\vert  \U\vert^pdx\right)^{\frac{1}{p}}.$$ 
Recalling that $\frac{R}{2}\leq \rho<r<R$ we finally get
\begin{equation}\label{estimate_I_1}
|I_{1}| \leq  c\frac{R^{3\left( \frac{1}{2}-\frac{1}{p}\right)}}{r-\rho}\left(\int_{B_r}\vert \vec{\nabla}\otimes \U\vert^2dx\right)^{\frac{1}{2}} \left(\int_{B_r}\vert  \U\vert^pdx\right)^{\frac{1}{p}}. 
\end{equation}
We study now the term $I_2$ in the identity (\ref{identity_ipp}).  By the Cauchy-Schwarz inequality we write
\begin{eqnarray*}
\vert I_2 \vert &=&\left|  \sum_{i,j=1}^{3} \int_{B_r} \partial_j U_i \partial_j (W_R)_i dx\right| \leq   c \left(\int_{B_r}\vert \vec{\nabla}\otimes \U\vert^2dx\right)^{\frac{1}{2}}\left(\int_{B_r}\vert \vec{\nabla}\otimes \W_R\vert^2dx\right)^{\frac{1}{2}}\\
&\leq & c \left(\int_{B_r}\vert \vec{\nabla}\otimes \U\vert^2dx\right)^{\frac{1}{2}} r^{3\left(\frac{1}{2}-\frac{1}{p}\right)}   \left(\int_{B_r}\vert \vec{\nabla}\otimes \W_R\vert^pdx\right)^{\frac{1}{p}}.
\end{eqnarray*}
But, by estimate (\ref{estim-W}) we have $\ds{\left(\int_{B_r}\vert \vec{\nabla}\otimes \W_R\vert^pdx\right)^{\frac{1}{p}} \leq c \left(\int_{B_r}\vert \vec{\nabla}\varphi_R\cdot \U \vert^pdx\right)^{\frac{1}{p}}}$, and by the last estimate  we can write
$$ \vert I_2 \vert \leq  c \left(\int_{B_r}\vert \vec{\nabla}\otimes \U\vert^2dx\right)^{\frac{1}{2}} r^{3\left(\frac{1}{2}-\frac{1}{p}\right)}\left(\int_{B_r}\vert \vec{\nabla}\varphi_R\cdot \U \vert^pdx\right)^{\frac{1}{p}}.$$ 
Again, since $\frac{R}{2}\leq \rho<r<R$ we have
\begin{equation}\label{estmate_I_2}
\vert I_2 \vert \leq  c\frac{R^{3\left( \frac{1}{2}-\frac{1}{p}\right)}}{r-\rho}\left(\int_{B_r}\vert \vec{\nabla}\otimes \U\vert^2dx\right)^{\frac{1}{2}} \left(\int_{B_r}\vert  \U\vert^pdx\right)^{\frac{1}{p}}.
\end{equation}
With inequalities (\ref{estimate_I_1}) and (\ref{estmate_I_2}), the Lemma \ref{Lemme1} is proven.   \finpv  
In order to study the last quantity $I_{3}$ in (\ref{identity_ipp}) we will need the following lemma.
 \begin{Lemme}\label{Lemme2} Let $\frac{9}{2}<p<6$ and let $\U \in L^p\cap \dot{B}^{\frac{3}{p}-\frac{3}{2},\infty}_{\infty}(\Rt)$ be a weak solution of the stationary Navier-Stokes equations  (\ref{N-S-f-nulle-intro}). Then we have the following estimate:  
 \begin{equation}\label{estim-I_3}	
 \ds{\vert I_3 \vert \leq   c\frac{R}{r-\rho}   \Vert \U \Vert_{\dot{B}^{\frac{3}{p}-\frac{3}{2}, \infty}_{\infty}}    \left(\int_{B_r}\vert \vec{\nabla}\otimes \U\vert^2dx\right)^{\frac{1}{2}} \left(\int_{B_r}\vert  \U\vert^pdx\right)^{\frac{1}{p}}} 
 \end{equation} where $c>0$ is always a constant which does not depend of $R,r,\rho$ and  $\U$. 
\end{Lemme}
This lemma is technical and we postpone to the appendix the details of its proof.\\

Thus, by equation (\ref{identity_ipp}) and with the inequalities of Lemmas \ref{Lemme1} and \ref{Lemme2}, we can write 
\begin{eqnarray} \label{estim-aux3} \nonumber
\int_{B_r}\varphi_R \vert \vec{\nabla}\otimes \U\vert^2dx & \leq c & \left( \frac{R^{3\left( \frac{1}{2}-\frac{1}{p}\right)}}{r-\rho} + \frac{R}{r-\rho}  \Vert \U \Vert_{\dot{B}^{\frac{3}{p}-\frac{3}{2}, \infty}_{\infty}} \right) \left(\int_{B_r}\vert \vec{\nabla}\otimes \U\vert^2dx\right)^{\frac{1}{2}} 	  \\
&  & \times  \left(\int_{B_r}\vert  \U\vert^pdx\right)^{\frac{1}{p}}.
\end{eqnarray}
Moreover,  for the first term in the right side we have 
$$ c \left( \frac{R^{3\left( \frac{1}{2}-\frac{1}{p}\right)}}{r-\rho} + \frac{R}{r-\rho}  \Vert \U \Vert_{\dot{B}^{\frac{3}{p}-\frac{3}{2}, \infty}_{\infty}} \right) \leq c  \  \left( \frac{R^{3\left( \frac{1}{2}-\frac{1}{p}\right)}+R}{r-\rho}  \right)  \left(  1+\Vert \U \Vert_{\dot{B}^{\frac{3}{p}-\frac{3}{2}, \infty}_{\infty}}\right),$$ and we set now the constant 
\begin{equation}\label{c(R)}
\mathfrak{C}(\U,R)=  c  \left(R^{3\left( \frac{1}{2}-\frac{1}{p}\right)}+R\right) \left(  1+\Vert \U \Vert_{\dot{B}^{\frac{3}{p}-\frac{3}{2}, \infty}_{\infty}}\right)>0, 
\end{equation} and thus we write 
$$ c \left( \frac{R^{3\left( \frac{1}{2}-\frac{1}{p}\right)}}{r-\rho} + \frac{R}{r-\rho}  \Vert \U \Vert_{\dot{B}^{\frac{3}{p}-\frac{3}{2}, \infty}_{\infty}} \right) \leq \frac{\mathfrak{C}(\U,R)}{r-\rho}.$$
With these estimates we get back to inequality (\ref{estim-aux3}) and we have the following estimate:

$$ \int_{B_r}\varphi_R\vert \vec{\nabla}\otimes \U\vert^2dx\leq \frac{\mathfrak{C}(\U,R)}{(r-\rho)} \left(\int_{B_r}\vert \vec{\nabla}\otimes \U\vert^2dx\right)^{\frac{1}{2}} \Vert \U \Vert_{L^p}.$$ 
On the other hand, as  $\varphi_R(x)=1$ if $\vert x \vert <\rho$, we have $\ds{\int_{B_{\rho}}\vert \vec{\nabla}\otimes \U \vert^2 dx \leq \int_{B_r}\varphi_R\vert \vec{\nabla}\otimes \U\vert^2dx}$, and by the last estimate we can write
$$\int_{B_{\rho}}\vert \vec{\nabla}\otimes \U \vert^2 \leq \frac{\mathfrak{C}(\U,R)}{(r-\rho)} \left(\int_{B_r}\vert \vec{\nabla}\otimes \U\vert^2dx\right)^{\frac{1}{2}} \Vert \U \Vert_{L^p},$$
where, applying the Young inequalities (with $1=\frac{1}{2}+\frac{1}{2}$) in the term in the right side we obtain the following inequality
\begin{equation}\label{Young}
\int_{B_{\rho}}\vert \vec{\nabla}\otimes \U \vert^2 dx \leq  \frac{1}{4} \int_{B_r}\vert \vec{\nabla}\otimes \U\vert^2dx +4\frac{\mathfrak{C}^{2}(\U,R) \Vert \U \Vert^{2}_{L^p}}{(r-\rho)^2}.
\end{equation} 
With this inequality  at hand, we obtain the desired estimate (\ref{Cacciopoli}) as follows: for all  $k \in \mathbb{N}$ positive we set  $ \ds{\rho_k=\frac{R}{2^{\frac{1}{k}}}}$, and in  estimate  (\ref{Young}) we set  $\rho=\rho_k$ and $r=\rho_{k+1}$ (where $\frac{R}{2} \leq \rho_k < \rho_{k+1} <R$) and then  we write  
\begin{equation}\label{estim-aux4} 
\int_{B_{\rho_k}}\vert \vec{\nabla}\otimes \U \vert^2 \leq \frac{1}{4}\int_{B_{\rho_{k+1}}}\vert \vec{\nabla}\otimes \U\vert^2dx+ 4\frac{\mathfrak{C}^{2}(\U,R) \Vert \U \Vert^{2}_{L^p}}{(\rho_{k+1}-\rho_k)^2}.
	\end{equation}
Now, let us study  the second term in the right side. Since  $ \ds{\rho_k=\frac{R}{2^{\frac{1}{k}}}}$ then we have $\ds{(\rho_{k+1}-\rho_k)^2 =R^2 \left( \frac{1}{2^{\frac{1}{k+1}}} - \frac{1}{2^{\frac{1}{k}}} \right)^2}$. But, for $k\in \mathbb{N}$ positive we have  $\ds{\frac{1}{2^{\frac{1}{k+1}}} - \frac{1}{2^{\frac{1}{k}}} \geq c \frac{1}{k}}$, where  $c>0$ is a numerical constant which does not depend of $k$, and thus we have $\ds{ (\rho_{k+1}-\rho_k)^2 \geq c \frac{R^2}{k^2}}$, hence we  write  
	$$ 4\frac{\mathfrak{C}^{2}(\U,R) \Vert \U \Vert^{2}_{L^p}}{(\rho_{k+1}-\rho_k)^2} \leq  4c\,k^2 \frac{\mathfrak{C}^{2}(\U,R) \Vert \U \Vert^{2}_{L^p}}{R^2}.$$ Then, with this estimate and  getting back to inequality (\ref{estim-aux4}) we get the following recursive formula: 
	\begin{equation*}\label{recurrence}
 \int_{B_{\rho_k}}\vert \vec{\nabla}\otimes \U \vert^2 \leq \frac{1}{4}\int_{B_{\rho_{k+1}}}\vert \vec{\nabla}\otimes \U\vert^2dx +4c\,k^2 \frac{\mathfrak{C}^{2}(\U,R) \Vert \U \Vert^{2}_{L^p}}{R^2}.
	\end{equation*}
Now, iterating this recursive formula for  $k=1,\cdots, n$ and since  $\rho_{1}=\frac{R}{2}$ we get the following estimate
\begin{equation*}\label{estim-aux5}
\int_{B_{\frac{R}{2}}}\vert \vec{\nabla}\otimes \U \vert^2dx\leq \frac{1}{4^n}\int_{B_{\rho_{n+1}}}\vert \vec{\nabla}\otimes \U\vert^2dx +4c  \frac{\mathfrak{C}^{2}(\U,R) \Vert \U \Vert^{2}_{L^p}}{R^2} \left( \sum_{k=1}^{n} \frac{k^2}{4^{k-1}}\right).
\end{equation*}	
In this estimate, recall that $\rho_{n+1} < R$  and then we can write 
$$ \int_{B_{\frac{R}{2}}}\vert \vec{\nabla}\otimes \U \vert^2dx\leq \frac{1}{4^n}\int_{B_{R}}\vert \vec{\nabla}\otimes \U\vert^2dx +4c \, \frac{\mathfrak{C}^{2}(\U,R) \Vert \U \Vert^{2}_{L^p}}{R^2} \left( \sum_{k=1}^{n} \frac{k^2}{4^{k-1}}\right),$$	
and  taking  the limit when $n\longrightarrow +\infty$ and  since $\ds{\sum_{k=1}^{+\infty} \frac{k^2}{4^{k-1}} <+\infty}$ then we have 
\begin{equation}\label{estim-1}
\int_{B_{\frac{R}{2}}}\vert \vec{\nabla}\otimes \U \vert^2dx\leq c \frac{\mathfrak{C}^{2}(\U,R) \Vert \U \Vert^{2}_{L^p}}{R^2}.
\end{equation} Finally, in this inequality we study the term $\ds{ \frac{\mathfrak{C}^{2}(\U,R)}{R^2}}$. Recall that the quantity $\mathfrak{C}(\U,R)$ is defined in expression 
(\ref{c(R)}) and by this expression we have 
\begin{eqnarray*} 
\frac{\mathfrak{C}^{2}(\U,R)}{R^2} & \leq & c \frac{1}{R^2} \left(  R^{3\left( \frac{1}{2}-\frac{1}{p}\right)}+R\right)^{2}  \left(  1+\Vert \U \Vert_{\dot{B}^{\frac{3}{p}-\frac{3}{2}, \infty}_{\infty}}\right)^{2} \\
& \leq &  c \left( R^{1-\frac{6}{p}} +1 \right)  \left(  1+\Vert \U \Vert_{\dot{B}^{\frac{3}{p}-\frac{3}{2}, \infty}_{\infty}}\right)^2.
\end{eqnarray*}   Thus, we define now the constant $\ds{C(\U,R)=c  \left(R^{1-\frac{6}{p}}+1 \right)   \left(1+ \Vert \U \Vert_{\dot{B}^{\frac{3}{p}-\frac{3}{2},\infty}_{\infty}}\right)^{2}}$ and by estimate (\ref{estim-1}) we have  Cacciopoli type  estimate (\ref{Cacciopoli}). \\
\\	 
With the estimate (\ref{Cacciopoli}) we can prove now that $\U \in \dot{H}^1(\Rt)$. Indeed, by this estimate we can write 
$$ \int_{B_{\frac{R}{2}}}\vert \vec{\nabla}\otimes \U (x)\vert^2dx \leq C(\U, R)  \Vert \U\Vert^{2}_{L^p}.$$ But, since $\ds{C(\U,R)=c \left(R^{1-\frac{6}{p}}+1 \right)   \left(1+ \Vert \U \Vert_{\dot{B}^{\frac{3}{p}-\frac{3}{2},\infty}_{\infty}}\right)^{2}}$, and  since   $\frac{9}{2}<p<6$ then we have $-\frac{1}{3}<1-\frac{6}{p}<0$ and thus we can write  $$\ds{\lim_{R\longrightarrow +\infty} C(\U,R)= c\left(1+ \Vert \U \Vert_{\dot{B}^{\frac{3}{p}-\frac{3}{2},\infty}_{\infty}}\right)^{2} <+\infty}.$$ Now,   in estimate estimate  (\ref{Cacciopoli}) we take the limit when $R\longrightarrow +\infty$  and we get  $ \Vert \U \Vert^{2}_{\dot{H}^1}\leq c\left(1+ \Vert \U \Vert_{\dot{B}^{\frac{3}{p}-\frac{3}{2},\infty}_{\infty}}\right)^{2}  \Vert \U \Vert^{2}_{L^p}<+\infty$. Proposition \ref{Prop:Lp-Besov} is now proven. \finpv  
 \\
By  Proposition \ref{Prop:Lp-Besov} we have the information   $\U\in \dot{H}^1(\Rt)$ and now  we can prove the identity $\U=0$.  Recall that we also have the information $\U \in \dot{B}^{\frac{3}{p}-\frac{3}{2},\infty}_{\infty}(\Rt)$ and then if we set the parameter $\beta=\frac{3}{2}-\frac{3}{p}$ (where, as  $\frac{9}{2}<p<6$ then we have $\frac{5}{6}<\beta <1$) then by the improved Sobolev inequalities  (see  the article  \cite{GerardMeyerOru}) we can write  
\begin{equation}\label{sobolev-presice}
\Vert \U  \Vert_{L^q} \leq c \Vert  \U \Vert^{\theta}_{\dot{H}^{1}} \Vert \U \Vert^{1-\theta}_{\dot{B}^{-\beta, \infty}_{\infty}},
\end{equation} with $\theta=\frac{2}{q}$ and   $\beta=\frac{\theta}{1-\theta}$; and by these identities  we have the following relation  $q=\frac{2}{\beta}+2$, where, as $\frac{5}{6}<\beta <1$ then we have  $3<q<\frac{9}{2}$. \\
\\
Once we have $\U \in L^q(\Rt)$, with   $3< q<\frac{9}{2}$, by  point $1)$ of Theorem \ref{Theo:Lebesgue} we can write $\U=0$. This finish the proof of the second point  of Theorem \ref{Theo:Lebesgue} and this theorem is now proven. \finpv   
\end{enumerate}


\section{The Liouville problem in Morrey spaces} \label{Ths-Morrey}
In this section we study the \emph{Liouville problem} for the stationary Navier-Stokes equations (\ref{N-S-f-nulle-intro}) where the weak solution $\U \in L^{2}_{loc}(\Rt)$ belongs to  Morrey spaces.

\subsection{Proof of Theorem \ref{Theo:Seregin-generalise}}
Assume that $\U\in \dot{M}^{2,3}(\Rt)\cap \dot{M}^{2,q}(\Rt)$ with $3<q<+\infty$. We will prove the identity $\U=0$ and for this, first we need to prove that the solution $\U$ also belongs to the Lebesgue space $L^{\infty}(\Rt)$. \\
\\
Indeed, let us consider the stationary solution  $\U \in \dot{M}^{2,q}(\Rt)$ as the initial data of the  Cauchy problem for the non  stationary  Navier-Stokes equations:
\begin{equation}\label{N-S-aux2}
\partial_t \vu + (\vu \cdot \vec{\nabla})\vu - \Delta \vu +\vec{\nabla} p =0, \quad div(\vu)=0, \quad \vu(0,\cdot)=\U.
\end{equation} 
By Theorem  $8.2$ (page $166$) of the book \cite{PGLR1}, there exists a time $T_0>0$, and a function  $\vu\in \mathcal{C}([0,T_0[,\dot{M}^{2,q}(\Rt))$ which is a solution of the Cauchy problem  (\ref{N-S-aux2}) and which also verifies the estimate
\begin{equation}\label{estim:Linf}
\sup_{0<t<T_0} t^{\frac{3}{2 q}}\Vert \vu(t,\cdot)\Vert_{L^{\infty}}<+\infty. 
\end{equation}
Moreover, by Theorem $8.4$ (page $172$) of book the \cite{PGLR1}, for the values $3<q<+\infty$ we have the uniqueness of  this solution $ \vu\in \mathcal{C}([0,T_0[,\dot{M}^{2,q}(\Rt))$.  But, since $\U \in \dot{M}^{2,q}(\Rt)$ is a stationary function  then we have $\U \in \mathcal{C}([0,T_0[,\dot{M}^{2,q}(\Rt))$ and since $\U$ is a solution of the stationary Navier-Stokes equations (\ref{N-S-f-nulle-intro}) then this function  is also a solution for the Cauchy  problem (\ref{N-S-aux2}) (since we have $\partial_t \U =0$) and thus, by uniqueness of solution $\vu$, we have the identity  $\vu=\U$. \\
\\
Thus, by estimate  (\ref{estim:Linf}) we can write  
\begin{equation}\label{estim-U-Linf}
\left(\frac{T_0}{2}\right)^{\frac{3}{2 q}} \Vert \U \Vert_{L^{\infty}}\leq \sup_{0<t<T_0} t^{\frac{1}{2}}\Vert \U\Vert_{L^{\infty}}<+\infty, 
\end{equation} and we get  $\U \in L^{\infty}(\Rt)$. \\
\\
Once we have the information $\U \in L^{\infty}(\Rt)$, we will use the additional information  $\U \in \dot{M}^{2,3}(\Rt)$ in order to prove $\U =0$. Let us start by proving the following proposition: 
\begin{Proposition}\label{Prop:morrey-sobolev} Let $\U \in L^{\infty} \cap  \dot{M}^{2,3}(\Rt)$ be a solution of stationary Navier-Stokes equations  (\ref{N-S-f-nulle-intro}). Then  $\U \in \dot{H}^{1}(\Rt)$ and we have $\Vert \U \Vert_{\dot{H}^{1}}\leq c\Vert \U \Vert^{\frac{1}{2}}_{L^{\infty}}\Vert \U \Vert_{\dot{M}^{2,3}}$. 
\end{Proposition}
\pv  Let $R>1$ and $B_R=\{ x \in \Rt: \vert x \vert <R\}$. We will prove the following estimate 
\begin{equation}\label{estim-aux7}
\int_{B_\frac{R}{2}}\vert \vec{\nabla} \otimes \U\vert^{2}dx \leq c \left( \frac{1}{R^{\frac{1}{3}}} \Vert \U \Vert^{\frac{2}{3}}_{L^{\infty}} \Vert \U \Vert^{\frac{4}{3}}_{\dot{M}^{2,3}} + \Vert \U \Vert_{L^{\infty}} \Vert \U \Vert^{2}_{\dot{M}^{2,3}} \right).
\end{equation}
For this,  following some ideas of the articles \cite{Ser2015} and \cite{Ser2016}, the first thing to do is to define  the following cut-off function: for a fixed $R>1$, we define the function $\phi_R \in \mathcal{C}^{\infty}_{0}(\Rt)$ such that $0\leq \phi_R \leq 1$, $\phi_R(x)=1$ if $\vert  x \vert < \frac{R}{2}$, $\phi_R(x)=0$ if $\vert x \vert >R$  and moreover  this function verifies  $\ds{ \Vert \vec{\nabla} \phi_R \Vert_{L^{\infty}}\leq \frac{c}{R}}$ and  $\ds{ \Vert \Delta \phi_R \Vert_{L^{\infty}}\leq \frac{c}{R^2}}$, where $c>0$ is a constant which does not depend of $R>1$. \\
\\
With this function  $\phi_R$ and the stationary  solution $\U $ we consider now the function  $\phi_R \U$ and we write 
\begin{equation}\label{iden-ipp-1}
\int_{B_R} \left(  -\Delta \U +(\U \cdot \vec{\nabla})\U+\vec{\nabla}P\right)\cdot (\phi_R \U )dx=0,
\end{equation} 
Now, we must study this identity and for this  we need first the following technical lemma:
\begin{Lemme}\label{Lemme:morrey-L-inf-interp} Let $\U \in L^{\infty} \cap  \dot{M}^{2,3}(\Rt)$. Then  we have $\Vert \U \Vert_{\dot{M}^{3,\frac{9}{2}}}\leq c \Vert \U \Vert^{\frac{1}{3}}_{L^{\infty}} \Vert U \Vert^{\frac{2}{3}}_{\dot{M}^{2,3}}$. \end{Lemme}
\pv 
Let  $x_0\in \Rt$ and $r>0$. Let  the ball $B(x_0,R)\subset \Rt$, we have  
$$ \left( \int_{B(x_0,r)} \vert \U \vert^3 dx\right)^{\frac{1}{3}} \leq c \left(\left(\int_{B(x_0,r)} \vert \U \vert^2 dx \right)^{\frac{1}{2}} \right)^{\frac{2}{3}} \Vert \U \Vert^{\frac{1}{3}}_{L^{\infty}}, $$ and multiplying by $r^{-\frac{1}{3}}$ in both sides of this estimate we get 
\begin{eqnarray*}
r^{-\frac{1}{3}}  \left( \int_{B(x_0,r)} \vert \U \vert^3 dx\right)^{\frac{1}{3}} & \leq &c \,   r^{-\frac{1}{3}}  \left( \left(\int_{B(x_0,r)} \vert \U \vert^2 dx \right)^{\frac{1}{2}} \right)^{\frac{2}{3}} \Vert \U \Vert^{\frac{1}{3}}_{L^{\infty}} \\
& \leq & c \left(r^{-\frac{1}{2}}  \left(\int_{B(x_0,r)} \vert \U \vert^2 dx \right)^{\frac{1}{2}} \right)^{\frac{2}{3}} \Vert \U \Vert^{\frac{1}{3}}_{L^{\infty}}.
\end{eqnarray*}
Now, if  in the first estimate in the  left side we write $\ds{r^{-\frac{1}{3}}= r^{\frac{3}{\frac{9}{2}}-\frac{3}{3} }}$ and moreover,  if  in the last estimate to the right side we write  $\ds{r^{-\frac{1}{2}}=r^{\frac{3}{3} -\frac{3}{2}}}$, then we have
$$  r^{\frac{3}{\frac{9}{2}}-\frac{3}{3} }   \left( \int_{B(x_0,r)} \vert \U \vert^3 dx\right)^{\frac{1}{3}} \leq c  \left(r^{\frac{3}{3} -\frac{3}{2}}  \left(\int_{B(x_0,r)} \vert \U \vert^2 dx \right)^{\frac{1}{2}} \right)^{\frac{2}{3}} \Vert \U \Vert^{\frac{1}{3}}_{L^{\infty}},$$  and thus we can write 
$$ \sup_{x_0 \in \Rt, r>0} \left(  r^{\frac{3}{\frac{9}{2}}-\frac{3}{3} }   \left( \int_{B(x_0,r)} \vert \U \vert^3 dx\right)^{\frac{1}{3}} \right) \leq c \left( \sup_{x_0 \in \Rt, r>0}  \left(  r^{\frac{3}{3} -\frac{3}{2}}  \left(\int_{B(x_0,r)} \vert \U \vert^2 dx \right)^{\frac{1}{2}} \right) \right)^{\frac{2}{3}} \Vert \U \Vert^{\frac{1}{3}}_{L^{\infty}}.$$ Finally, by definition of quantities  $\Vert \U \Vert_{\dot{M}^{3,\frac{9}{2}}}$ and $\Vert \U \Vert_{\dot{M}^{2,3}}$ given in formula (\ref{def-morrey}) we can write $\Vert \U \Vert_{\dot{M}^{3,\frac{9}{2}}}\leq c  \Vert \U \Vert^{\frac{1}{3}}_{L^{\infty}} \Vert U \Vert^{\frac{2}{3}}_{\dot{M}^{2,3}}$. \finpv
\\
Once we have the information $\U \in \dot{M}^{3,\frac{9}{2}}(\Rt)$ we get back to study the identity (\ref{iden-ipp-1}). \\
\\
Remark first that since $\U \in \dot{M}^{3,\frac{9}{2}}(\Rt)$  then we have  $\U \in L^{3}_{loc}(\Rt)$  and thus, by Theorem X.1.1 of the book \cite{Galdi} (page $658$),  we have  $\U \in \mathcal{C}^{\infty}(\Rt)$ and $P\in \mathcal{C}^{\infty}(\Rt)$ and thus  all the terms in (\ref{iden-ipp-1}) are well-defined and they are smooth enough.\\
\\
Then, we can integrate by parts each term in the identity (\ref{iden-ipp-1}): for the first term $\ds{\int_{B_R} \left(  -\Delta \U \right) \cdot ( \phi_R \U)dx }$,   following the same computations  in equation (\ref{ipp1}) (with the function $\phi_R$ in instead of the function $\theta_R$)   we have
$$ \int_{B_R} \left(  -\Delta \U \right) \cdot ( \phi_R \U)dx= -  \int_{B_R} \Delta \phi_R  \left( \frac{\vert \U\vert^2}{2}  \right)dx +   \int_{B_R}\phi_R \vert \vec{\nabla}\otimes \U \vert^{2} dx. $$
For the second term in identity (\ref{iden-ipp-1}): $\ds{\int_{B_R} \left( (\U \cdot \vec{\nabla})\U\right)\cdot (\phi_R \U )dx}$, always following the same computations in equations (\ref{ipp2}) and (\ref{ipp3}) we can write  
$$ \int_{B_R} \left( (\U \cdot \vec{\nabla})\U\right)\cdot (\phi_R \U )dx= - \int_{B_R} \vec{\nabla} \phi_R \cdot \left( \frac{\vert \U\vert^2}{2} \U\right)dx.  $$
Finally, for the third term in identity (\ref{iden-ipp-1}): $\ds{\int_{B_R} \left(  \vec{\nabla}P\right)\cdot (\phi_R \U )dx}$, following again the same computations as in equation (\ref{ipp4}) we have 
$$ \int_{B_R} \left(  \vec{\nabla}P\right)\cdot (\phi_R \U )dx = - \int_{B_R}\vec{\nabla}\phi_R \cdot(P \U) dx. $$
With these identities, we get back to the identity (\ref{iden-ipp-1}) and we write 
$$   - \int_{B_R} \Delta \phi_R  \left( \frac{\vert \U\vert^2}{2}  \right)dx +   \int_{B_R}\phi_R \vert \vec{\nabla}\otimes \U \vert^{2} dx - \int_{B_R} \vec{\nabla} \phi_R \cdot \left( \frac{\vert \U\vert^2}{2} \U\right)dx  - \int_{B_R}\vec{\nabla}\phi_R \cdot(P \U) dx=0,$$ hence we have
 
\begin{eqnarray}\label{iden1-aux} \nonumber
\int_{B_R} \phi_R \vert \vec{\nabla}\otimes \U \vert^2 dx &= &\int_{B_R} \Delta \phi_R \frac{\vert \U \vert^2}{2} dx + \int_{B_R} \vec{\nabla} \phi_R \cdot \left( \left( \frac{\vert \U \vert^2}{2}+P \right) \U \right) dx\\
&=&  I_1(R)+I_2(R),
\end{eqnarray}
and  we study now the terms   $I_1(R)$ and $I_2(R)$.\\
\\
\\
For the first term  $I_1(R)$, as  $\ds{ \Vert \Delta \phi_R \Vert_{L^{\infty}}\leq \frac{c}{R^2}}$ we have
$$  \vert I_1(R) \vert \leq \int_{B_R} \vert \Delta \phi_R \vert \frac{\vert \U\vert^2}{2} dx \leq \frac{c}{R^2}\int_{B_R}\vert \U \vert^2 dx,$$ and in the last term in the right side  we can write 
$$ \frac{c}{R^2}\int_{B_R}\vert \U \vert^2 dx \leq \frac{c}{R^2}  \left( R^{6\left(\frac{1}{2}-\frac{1}{3} \right)}   \left( \int_{B_R}\vert \U \vert^3dx\right)^{\frac{2}{3}} \right)\leq \frac{c}{R}\left( \int_{B_R}\vert \U \vert^3dx\right)^{\frac{2}{3}}.$$ 
But,  since  $\U \in \dot{M}^{3,\frac{9}{2}}(\Rt)$ 
then by expression (\ref{def-morrey}) we have  $$\left( \int_{B_R}\vert \U \vert^3dx\right)^{\frac{2}{3}} \leq R^{6\left( \frac{1}{3}-\frac{2}{9}\right)}  \Vert \U \Vert^{2}_{\dot{M}^{3,\frac{9}{2}}},$$ and thus we get 
$$ \frac{c}{R} \left( \int_{B_R}\vert \U \vert^3dx\right)^{\frac{2}{3}} \leq c\,\frac{R^{6\left( \frac{1}{3}-\frac{2}{9}\right)}}{R}   \Vert \U \Vert^{2}_{\dot{M}^{3,\frac{9}{2}}}= \frac{c}{R^{\frac{1}{3}}} \Vert \U \Vert^{2}_{\dot{M}^{3,\frac{9}{2}}}. $$
Thus, by these estimates we finally  get 
\begin{equation}\label{estim_I_1_Prop_Th_2}
\vert I_1(R) \vert   \leq \frac{c}{R^{\frac{1}{3}}}   \Vert \U \Vert^{2}_{\dot{M}^{3,\frac{9}{2}}}. 
\end{equation} 
For the second term $I_2(R)$ in (\ref{iden1-aux}),  since  $\ds{ \Vert \vec{\nabla} \phi_R \Vert_{L^{\infty}}\leq \frac{c}{R}}$ then   we  can write
\begin{eqnarray}\label{estim-aux-morrey} \nonumber
\vert I_2(R) \vert & \leq &\int_{B_R} \vert \vec{\nabla} \phi_R  \vert \left\vert  \left( \frac{\vert \U \vert^2}{2}+P \right) \U \right\vert dx  \leq  \frac{c}{R}\int_{B_R}\vert \U \vert^3 dx+ \frac{c}{R} \int_{B_R} \vert P \vert \vert \U \vert dx \\
&\leq& (I_2)_a+(I_2)_b,
\end{eqnarray}  and we still need to study the terms  $(I_2)_a$ and $(I_2)_b$ above.\\
\\
In order to study the  term $(I_2)_a$, recall first that  $\U \in \dot{M}^{3,\frac{9}{2}}(\Rt)$ and  by expression (\ref{def-morrey})  we can  write  $\ds{ \int_{B_R}\vert \U \vert^3dx\leq  R^{9\left(\frac{1}{3}-\frac{2}{9}\right)}   \Vert \U\Vert^{3}_{\dot{M}^{3,\frac{9}{2}}}}$.   Thus we get  
\begin{equation}\label{estim-I2a}
(I_2)_a\leq\ds{\frac{c}{R}\int_{B_R}\vert \U \vert^3 dx \leq  \frac{c}{R}  \left( R^{9\left(\frac{1}{3}-\frac{2}{9}\right)}   \Vert \U\Vert^{3}_{\dot{M}^{3,\frac{9}{2}}}\right)  \leq c \Vert \U\Vert^{3}_{\dot{M}^{3,\frac{9}{2}}}}.
\end{equation}
For the term  $(I_2)_b$,  applying the H\"older inequalities (with $1=\frac{2}{3}+ \frac{1}{3}$) we can write 
\begin{equation}\label{estim-aux6}
 \vert (I_2)_b \vert \leq  \frac{c}{R} \int_{B_R} \vert P \vert \vert \U \vert dx \leq \frac{c}{R}  \left(\int_{B_R}\vert P \vert^{\frac{3}{2}}dx \right)^{\frac{2}{3}} \left(\int_{B_R}\vert \U \vert^3 dx \right)^{\frac{1}{3}},
\end{equation} and we study now the two last terms in the right side. \\
\\
In order to estimate the  term $\ds{\left(\int_{B_R}\vert P \vert^{\frac{3}{2}}dx \right)^{\frac{2}{3}}}$ in the inequality  above  we need  the following technical lemma.
\begin{Lemme}\label{Lemme:pression} Let $(\U,P) \in L^{2}_{loc}(\Rt)\times \mathcal{D}^{'}(\Rt)$ be a solution of the stationary Navier-Stokes equations (\ref{N-S-f-nulle-intro}). If $\U \in  \dot{M}^{p,q}(\Rt)$ with $p\geq 2$ and $q\geq 3$ then we have $P \in \dot{M}^{\frac{p}{2},\frac{q}{2}}(\Rt)$ and $\Vert P \Vert_{\dot{M}^{\frac{p}{2},\frac{q}{2}}} \leq c \Vert \U \Vert^{2}_{\dot{M}^{p,q}}$.  
\end{Lemme} 
\pv By equation (\ref{pression}) we write the pressure $P$ as $\ds{P= \sum_{i,j=1}^{3}\mathcal{R}_i\mathcal{R}_j (U_i U_j)}$, where recall that  $\mathcal{R}_i=\frac{\partial_i}{\sqrt{-\Delta}}$
denotes the i-th Riesz transform. Then,  by continuity of the operator  $\mathcal{R}_i\mathcal{R}_j$ on Morrey spaces $\dot{M}^{p,q}(\Rt)$ for the values $p\geq 2$ and $q \geq 3$ (see the book  \cite{PGLR1}, page $171$) and applying the H\"older inequalities we get the following estimate
$$ \Vert P \Vert_{\dot{M}^{\frac{p}{2},\frac{q}{2}}}\leq c \sum_{i,j=1}^{3}\Vert \mathcal{R}_i \mathcal{R}_j(U_i U_j)\Vert_{\dot{M}^{\frac{p}{2},\frac{q}{2}}}\leq c \Vert \U \otimes \U \Vert_{\dot{M}^{\frac{p}{2},\frac{q}{2}}}\leq c \Vert \U \Vert^{2}_{\dot{M}^{p,q}}.$$ \finpv
Thus, since  $\U \in \dot{M}^{3,\frac{9}{2}}(\Rt)$ then by this lemma we have $P\in \dot{M}^{\frac{3}{2},\frac{9}{4}}(\Rt)$ and using the definition of the Morrey spaces given in (\ref{def-morrey}) we can write 
\begin{equation}\label{estim1}
\left(\int_{B_R}\vert P \vert^{\frac{3}{2}}dx \right)^{\frac{2}{3}} \leq   R^{3\left(\frac{2}{3}-\frac{4}{9} \right)}   \Vert P \Vert_{\dot{M}^{\frac{3}{2},\frac{9}{4}}}.
\end{equation}  For the term $\ds{\left(\int_{B_R}\vert \U \vert^3 dx \right)^{\frac{1}{3}}}$ in inequality (\ref{estim-aux6}), since $\U \in \dot{M}^{3,\frac{9}{2}}(\Rt)$ always by  expression (\ref{def-morrey})  we can write 
\begin{equation}\label{estim2}
\left(\int_{B_R}\vert \U \vert^3 dx \right)^{\frac{1}{3}} \leq   R^{3\left(\frac{1}{3}-\frac{2}{9} \right)}  \Vert \U \Vert_{\dot{M}^{3,\frac{9}{2}}}. 
\end{equation} Thus, with  estimates (\ref{estim1}) and (\ref{estim2}) we get back  to the inequality (\ref{estim-aux6}) and moreover, since by Lemma \ref{Lemme:pression} we have $\Vert P \Vert_{\dot{M}^{\frac{p}{2},\frac{q}{2}}} \leq c \Vert \U \Vert^{2}_{\dot{M}^{p,q}}$ then we obtain 

\begin{eqnarray}\label{estim_I2b} \nonumber 
\vert (I_2)_b \vert  &\leq &\frac{c}{R}   \left( R^{3\left(\frac{2}{3}-\frac{4}{9} \right)}   \Vert P \Vert_{\dot{M}^{\frac{3}{2},\frac{9}{4}}}  \right)\left(   R^{3\left(\frac{1}{3}-\frac{2}{9} \right)}  \Vert \U \Vert_{\dot{M}^{3,\frac{9}{2}}}  \right)\\ 
&\leq& c \Vert P \Vert_{\dot{M}^{\frac{3}{2},\frac{9}{4}}} \Vert \U \Vert_{\dot{M}^{3,\frac{9}{2}}}\leq c \Vert \U \Vert^{3}_{\dot{M}^{3,\frac{9}{2}}}.
\end{eqnarray}
Now, with estimates (\ref{estim-I2a}) and  (\ref{estim_I2b})  at hand, we get back to inequality  (\ref{estim-aux-morrey}) and we can  write
\begin{equation}\label{estim_term_I_2_Th_2}
\vert I_2(R)\vert\leq c \Vert \U \Vert^{3}_{\dot{M}^{3,\frac{9}{2}}}.
\end{equation} 
Once we have estimates (\ref{estim_I_1_Prop_Th_2}) and (\ref{estim_term_I_2_Th_2}), getting back to identity  (\ref{iden1-aux})  we have 
$$ \int_{B_{R}}\phi_R \vert \vec{\nabla}\otimes \U\vert^2dx \leq \frac{c}{R^{\frac{1}{3}}}  \Vert \U \Vert^{2}_{\dot{M}^{3,\frac{9}{2}}} +c \Vert \U \Vert^{3}_{\dot{M}^{3,\frac{9}{2}}}.$$ But, recall that   $\phi_R(x)=1$ if $\vert x \vert <\frac{R}{2}$ and then we have $\ds{\int_{B_{\frac{R}{2}}}\vert \vec{\nabla}\otimes \U\vert^2dx \leq \int_{B_{R}}\phi_R \vert \vec{\nabla}\otimes \U\vert^2dx }$  and thus we get the following estimate:
$$ \int_{B_{\frac{R}{2}}}\vert \vec{\nabla}\otimes \U\vert^2dx \leq  \frac{c}{R^{\frac{1}{3}}}  \Vert \U \Vert^{2}_{\dot{M}^{3,\frac{9}{2}}} +c \Vert \U \Vert^{3}_{\dot{M}^{3,\frac{9}{2}}}.$$
Moreover, recall that by Lemma \ref{Lemme:morrey-L-inf-interp} we have the estimate $\Vert \U \Vert_{\dot{M}^{3,\frac{9}{2}}}\leq c \Vert \U \Vert^{\frac{1}{3}}_{L^{\infty}} \Vert U \Vert^{\frac{2}{3}}_{\dot{M}^{2,3}}$, and thus we  finally obtain  the inequality (\ref{estim-aux7}). \\
\\
In order to finish the proof of Proposition \ref{Prop:morrey-sobolev}, in inequality (\ref{estim-aux7})  we take the limit  $R \longrightarrow+\infty$  and we  get  $\Vert \U \Vert_{\dot{H}^{1}}\leq c\Vert \U \Vert^{\frac{1}{2}}_{L^{\infty}}\Vert \U \Vert_{\dot{M}^{2,3}}$. \finpv
\\
\textbf{End of  the proof of Theorem \ref{Theo:Seregin-generalise}}. \\
\\
Now we have all the tools to prove the identity $\U=0$. First, recall that $\dot{M}^{2,3}(\Rt)$ is a  homogeneous Banach space  of degree $-1$ and then we have $\dot{M}^{2,3}(\Rt) \subset \dot{B}^{-1,\infty}_{\infty}(\Rt)$ (see the Chapter $4$ of the book \cite{PGLR2}). Thus, since $\U \in \dot{M}^{2,3}(\Rt)$ then
we have $\U \in \dot{B}^{-1,\infty}_{\infty}(\Rt)$. Moreover, by Proposition \ref{Prop:morrey-sobolev} we also have  $\U \in \dot{H}^{1}(\Rt)$ and then  by the improved Sobolev inequalities (\ref{sobolev-presice}) (with the parameters $\beta=1$, $\theta=\frac{1}{2}$ and $q=4$) we have $\U \in L^4(\Rt)$. Then, 
by point $1)$ of Theorem  \ref{Theo:Lebesgue} we can write $\U =0$ and Theorem  \ref{Theo:Seregin-generalise} is now proven. \finpv

\subsection{Proof of Theorem \ref{Theo:Morrey-2}}\label{Sec:Theo:Morrey-2} 
Assume here that the solution $\U \in L^{2}_{loc}(\Rt)$ of stationary Navier-Stokes equations (\ref{N-S-f-nulle-intro}) verifies $ \U\in \overline{M}^{p,3}(\Rt)$ with $2<p\leq 3$, where the space $ \overline{M}^{p,3}(\Rt)$ is given in Definition \ref{Def:sous-espace-morrey}. In order to prove the identity $\U=0$ we will follow some ideas of the proof of Theorem   \ref{Theo:Seregin-generalise} and the first thing to do is to prove that with this hypothesis on the solution $\U$ we have $\U \in L^{\infty}(\Rt)$. \\
\\
Indeed, we consider  the stationary solution $\U \in \overline{M}^{p,3}(\Rt)$  as the initial data of the  Cauchy problem for the non stationary Navier-Stokes equations (\ref{N-S-aux2}). Then, always by Theorem  $8.2$ of the book \cite{PGLR1}, there exists a function  $\vu \in \mathcal{C}([0,T_0[,\overline{M}^{p,3}(\Rt) )$ which is  a  solution of problem (\ref{N-S-aux2}). Moreover, this solution $\vu$ verifies the estimate:
\begin{equation}\label{estim} 
\sup_{0<t<T_0} t^{\frac{1}{2}}\Vert \vu(t,\cdot)\Vert_{L^{\infty}}<+\infty. 
\end{equation} 
On the other hand, recall that the stationary solution verifies $\U \in \mathcal{C}([0,T_0[,\overline{M}^{p,3}(\Rt) )$ and this function is also a solution of problem (\ref{N-S-aux2}) (always since $\partial_t \U=0$). But, for the values $2<p\leq 3$ by Theorem  $8.4$ of book \cite{PGLR1} we have the uniqueness of solution $\vu$ and thus we have the identity $\vu =\U$. By this identity we have that the function $\U$ verifies the estimate (\ref{estim}) hence, writing the same estimate as in equation (\ref{estim-U-Linf}), we get   $\U\in L^{\infty}(\Rt)$.  Remark here that Theorem $8.4$ assures  the uniqueness of solution $\vu$ in the space $\mathcal{C}([0,T_0[,\overline{M}^{p,3}(\Rt))$ and not in the more general setting of the space $\mathcal{C}([0,T_0[,\dot{M}^{p,3}(\Rt))$. For this reason we consider in Theorem \ref{Theo:Morrey-2} the functional space $\overline{M}^{p,3}(\Rt)$. \\
\\
We have now the information  $\U \in \overline{M}^{p,3} \cap L^{\infty}(\Rt)$ which will allows us to prove the  identity  $\U=0$. Indeed, recall that    $\overline{M}^{p,3} \subset \dot{M}^{2,3}(\Rt)$, hence we have $\U \in \dot{M}^{2,3}(\Rt) \cap L^{\infty}(\Rt)$ and by Proposition \ref{Prop:morrey-sobolev} we get $\U \in \dot{H}^{1}(\Rt)$. On the other hand,  since  $\dot{M}^{2,3}(\Rt) \subset \dot{B}^{-1,\infty}_{\infty}(\Rt)$ then the solution $\U \in \dot{M}^{2,3}(\Rt)$ verifies $\U \in \dot{B}^{-1,\infty}_{\infty}(\Rt)$  and the proof of the identity $\U=0$ follows the same lines given above at the end of the proof of Theorem \ref{Theo:Seregin-generalise}.  \finpv
\section{Appendix: Proof of Lemma \ref{Lemme2} page \pageref{Lemme1}}\label{appendix}
We prove here the estimate (\ref{estim-I_3}), where recall that the term $I_3$ (defined in the identity (\ref{identity_ipp}))  is given by  
\begin{equation}\label{I_3}
	I_3= -\int_{B_r}\left(U_1\partial_1\U+U_2\partial_2\U+U_3\partial_3\U\right) \cdot (\varphi_R \U-\W_R)dx.
	\end{equation}  
In order to study the term in right side above remark the  $\U$ can be written as   
	\begin{equation}\label{prim-U}
	\U=\vec{\nabla}\wedge \V
	\end{equation} where the vector field $\V$ is given by the following expression: 
	\begin{equation}\label{def-V}
	\V= \frac{1}{-\Delta} \left( \vec{\nabla}\wedge \U\right).
	\end{equation} Indeed, since $div(\U)=0$ then we have the following identities
	$$ \vec{\nabla}\wedge \V=\vec{\nabla}\wedge \left( \frac{1}{-\Delta} \left( \vec{\nabla}\wedge \U\right)  \right) =\frac{1}{-\Delta}  \left( \vec{\nabla} \wedge (\vec{\nabla} \wedge \U)\right) = \frac{1}{-\Delta} \left( \vec{\nabla} (div(\U))\right) - \frac{1}{-\Delta} \left( \Delta U\right) =\U.$$ 
But, in order  to carry out the estimates which we will need later, in equation (\ref{prim-U}) we will consider a little variant of function $\V$  above and we set now the function   $\V^{*}=\V-\V(0)$.  Remark that we have the identity  $\vec{\nabla}\wedge \V =\vec{\nabla}\wedge \V^{*}$ (because $\V(0)\in \Rt$ is a constant vector) and then by equation (\ref{prim-U}) we can write $\U=\vec{\nabla}\wedge  \V^{*}$, \emph{i.e.}, we have the identities $U_i=\partial_j V^{*}_{k}-\partial_{k}V^{*}_{j}$, where it is worth noting here that we always consider the indices $i,j,k\in \{1,2,3\}$  given by the right-hand rule:  if $i=1$ then $j=2$ and $k=2$; if $i=2$ then $j=3$ and $k=1$ and so on. \\
\\
Now, getting back to the term in the right side in  expression  (\ref{I_3}), we substitute  $U_i$ by $\partial_j V^{*}_{k}-\partial_{k}V^{*}_{j}$  and we write
	\begin{eqnarray*}
		I_3 &=&-\int_{B_r}\left((\partial_2 V^{*}_3-\partial_3 V^{*}_2)\partial_1\U+(\partial_3 V^{*}_1-\partial_1 V^{*}_3)\partial_2\U+(\partial_1 V^{*}_2-\partial_2 V^{*}_1)\partial_3\U\right)\cdot(\varphi_R \U-\W_R)dx\\
		&=&-\int_{B_r} \sum_{i=1}^{3} \left((\partial_jV^{*}_{k}-\partial_k V^{*}_{j})\partial_i \U \right)\cdot (\varphi_R \U-\W_R)dx\\
		&=& -\int_{B_r}  \sum_{i=1}^{3} \left(( \partial_jV^{*}_{k}) (\partial_i \U)\cdot (\varphi_R \U-\W_R)-( \partial_k V^{*}_{j})(\partial_i \U)\cdot (\varphi_R \U-\W_R)\right)dx.
	\end{eqnarray*} Then,   integrating by parts in each term above we have
	\begin{eqnarray*}
	I_3&=& \int_{B_r} \sum_{i=1}^{3}\left( \underbrace{V^{*}_{k}(\partial_j \partial_i \U )\cdot (\varphi_R \U-\W_R)}_{(a)}   + V^{*}_{k}(\partial_i\U)\cdot   \partial_j(\varphi_R \U-\W_R) \right)dx \\
	& & +\int_{B_r} \sum_{i=1}^{3}\left(  - \underbrace{V^{*}_{j}(\partial_k\partial_i\U)\cdot(\varphi_R \U-\W_R)}_{(b)} -V^{*}_{j}(\partial_i\U)\cdot\partial_k(\varphi_R \U-\W_R) \right)dx,\\
	\end{eqnarray*}
and grouping the terms 	$(a)$ and $(b)$ we can write 
	\begin{eqnarray}\label{estim:I3} \nonumber
I_3	&=& \int_{B_r} \sum_{i=1}^{3} \left( \left( V^{*}_{k}(\partial_j \partial_i \U )-V^{*}_{j}(\partial_k\partial_i\U)\right)\cdot (\varphi_R \U-\W_R)\right)dx\\ \nonumber
	& &  + \int_{B_r} \sum_{i=1}^{3} \left(V^{*}_{k}(\partial_i\U)\cdot   \partial_j(\varphi_R \U-\W_R) -V^{*}_{j}(\partial_i\U)\cdot\partial_k(\varphi_R \U-\W_R)\right)dx\\
	& =& (I_3)_a+(I_3)_b,
	\end{eqnarray}
where  we study now the terms  $(I_3)_a$ and $(I_3)_b$. \\
\\
For  the first term  $(I_3)_a$, recall that the  indices $i,j,k\in \{ 1,2,3\}$  are always given by the right-hand rule and then we have  $\ds{\sum_{i=1}^{3}  \left( V^{*}_{k}(\partial_j \partial_i \U )-V^{*}_{j}(\partial_k\partial_i\U)\right) =(0,0,0)}$ (just develop this sum to see that each term is canceled).  Thus we get  
	\begin{equation}\label{estim:Ia}
	(I_3)_a= \int_{B_r} \sum_{i=1}^{3} \left( \left( V^{*}_{k}(\partial_j \partial_i \U )-V^{*}_{j}(\partial_k\partial_i\U)\right)\cdot (\varphi_R \U-\W_R)\right)dx =0. 
	\end{equation}
	For the second term  $(I_3)_b$  we write 	
	\begin{eqnarray*}
		(I_3)_b &=& \int_{B_r}\sum_{i=1}^{3} V^{*}_{k}\partial_i \U\cdot  \left( (\partial_j \varphi_R) \U+ \underbrace{\varphi_R (\partial_j \U)}_{(c)} -\partial_j\W_R\right) dx   \\
		& & + \int_{B_r}\sum_{i=1}^{3} -V^{*}_{j}\partial_i\U \cdot \left( (\partial_k \varphi_R) \U+ \underbrace{\varphi_R (\partial_k \U)}_{(d)}-\partial_k \W_R\right) dx,
	\end{eqnarray*} 
and grouping now the terms $(c)$ and $(d)$ above we write 
\begin{eqnarray*}
	(I_3)_b &=&\int_{B_r} \sum_{i=1}^{3} \left( V^{*}_{k}\partial_i \U \cdot \left(\varphi_R (\partial_j\U) \right)-V^{*}_{j}\partial_i\U \cdot \left(\varphi_R(\partial_k\U)\right) \right) dx \\
	& & +\int_{B_r} \sum_{i=1}^{3} \left( V^{*}_{k}\partial_i\U\cdot (\partial_I\varphi_R \U-\partial_j\W_R) - V^{*}_{j}\partial_i \U \cdot(\partial_k \varphi_R \U-\partial_k \W_R)\right)dx.
\end{eqnarray*}
But, always since  the indices  $i,j,k \in \{1,2,3\}$ are given by the right-hand rule then we have $\ds{\sum_{i=1}^{3} \left( V^{*}_{k}\partial_i \U \cdot \left(\varphi_R (\partial_j\U) \right)-V^{*}_{j}\partial_i\U \cdot \left(\varphi_R(\partial_k\U)\right) \right)}=0$ (again,  develop this sum to see that each term is canceled)  and thus we get 
	\begin{equation}\label{estim:Ib}
	(I_3)_b=\int_{B_r} \sum_{i=1}^{3} \left( V^{*}_{k}\partial_i\U\cdot (\partial_i\varphi_R \U-\partial_j\W_R)- V^{*}_{j}\partial_i \U \cdot(\partial_k \varphi_R \U-\partial_k \W_R)\right)dx.
	\end{equation}
	With estimates (\ref{estim:Ia}) and (\ref{estim:Ib}),   we get back to term  $I_3$ given in identity (\ref{estim:I3}) and we can  write  
\begin{equation*}
I_3 = \int_{B_r} \sum_{i=1}^{3} \left( V^{*}_{k}\partial_i\U\cdot ( \underbrace{\partial_i \varphi_R \U}_{(e)} -\partial_j\W_R)- V^{*}_{j}\partial_i \U \cdot(\underbrace{\partial_k \varphi_R \U}_{(f)} -\partial_k \W_R)\right)dx,
\end{equation*}	where, 	grouping again the terms $(e)$ and $(f)$ above write  
	\begin{eqnarray*}
I_3&=& \int_{B_r} \sum_{i=1}^{3} \left( V^{*}_{k}\partial_i\U\cdot (\partial_i \varphi_R )\U- V^{*}_{j} \partial_i \U \cdot(\partial_k \varphi_R)\U \right)dx \\
& & + \int_{B_r} \sum_{i=1}^{3}\left( - V^{*}_{k}\partial_i\U\cdot (\partial_j \W_R)\U + V^{*}_{j}\partial_i\U\cdot \partial_k \W_R \right)dx, 
	\end{eqnarray*}  hence we have
\begin{eqnarray*}
\vert I_3 \vert  & \leq & 	\int_{B_r} \vert \V^{*}\vert \vert \vec{\nabla}\otimes \U\vert \vert \vec{\nabla}\varphi_R\otimes \U\vert dx + \int_{B_r} \vert \V^{*}\vert \vert \vec{\nabla}\otimes \U\vert \vert \vec{\nabla}\otimes \W_R\vert dx \\
&\leq & \int_{B_r} \vert \vec{\nabla}\otimes \U\vert  \left(  \vert \V^{*} \vert  \vert \vec{\nabla}\varphi_R\otimes \U\vert \right)  dx + \int_{B_r} \vert  \vec{\nabla}\otimes \U\vert   \left(  \vert \V^{*}\vert  \vert \vec{\nabla}\otimes \W_R\vert\right)dx. 
\end{eqnarray*}
 In both terms in the right side, applying first the  Cauchy-Schwarz inequality we write 
 \begin{eqnarray*}
\vert I_3 \vert  &\leq & \left( \int_{B_r} \vert \vec{\nabla} \otimes \U \vert^{2} dx \right)^{\frac{1}{2}}\left( \int_{B_r} \vert \V^{*}\vert^2 \vert  \vec{\nabla}\varphi_R\otimes \U\vert^2 dx \right)^{\frac{1}{2}}  \\
& &+ \left( \int_{B_r} \vert \vec{\nabla} \otimes \U \vert^{2} dx \right)^{\frac{1}{2}}\left( \int_{B_r}  \vert \V^{*}\vert^{2} \vert  \vec{\nabla}\otimes \W_R\vert^{2} dx \right)^{\frac{1}{2}},
\end{eqnarray*}
then, in each term in the right side we apply the  H\"older inequalities  (with $\frac{1}{2}=\frac{1}{q}+\frac{1}{p}$) and we have
$$\vert I_3 \vert \leq  \left(\int_{B_r}\vert \vec{\nabla}\otimes \U\vert^2dx\right)^{\frac{1}{2}}  \left(\int_{B_r}\vert \V^{*}\vert^{q}dx \right)^{\frac{1}{q}}   \left( \left( \int_{B_r} \vert \vec{\nabla}\varphi_R \otimes \U \vert^{p} dx \right)^{\frac{1}{p}} + \left(\int_{B_r} \vert \vec{\nabla} \otimes \W_R\vert^p dx \right)^{\frac{1}{p}} \right),$$ 
and we study now  the third term in the right side. Recall that by equation (\ref{control-test}) we have $\ds{ \Vert \vec{\nabla}\varphi_R \Vert_{L^{\infty}}\leq \frac{c}{r-\rho}}$ and then we can write 
$$ \left( \int_{B_r} \vert \vec{\nabla}\varphi_R \otimes \U \vert^{p} dx \right)^{\frac{1}{p}} \leq \frac{c}{r-\rho}  \left( \int_{B_r} \vert \U \vert^p dx\right)^{\frac{1}{p}}.$$
Moreover, recall that  by equation (\ref{estim-W}) we have 	$\ds{ \Vert \vec{\nabla}\otimes \W_R\Vert_{L^{p}(B_r)}\leq c \Vert \vec{\nabla}\varphi_R \cdot \U\Vert_{L^p(B_r)}}$ and then we have
$$ \left(\int_{B_r} \vert \vec{\nabla} \otimes \W_R\vert^p dx \right)^{\frac{1}{p}} \leq \frac{c}{r-\rho}  \left( \int_{B_r} \vert \U \vert^p dx\right)^{\frac{1}{p}}.$$ Thus, by these estimates we write  
 $$ \left( \left( \int_{B_r} \vert \vec{\nabla}\varphi_R \otimes \U \vert^{p} dx \right)^{\frac{1}{p}} + \left(\int_{B_r} \vert \vec{\nabla} \otimes \W_R\vert^p dx \right)^{\frac{1}{p}}\right) \leq \frac{c}{r-\rho} \left( \int_{B_r} \vert \U \vert^p dx\right)^{\frac{1}{p}}, $$ and then we get the following estimate 
\begin{equation}\label{iden-aux}
 \vert I_3 \vert \leq \frac{c}{r-\rho}  \left(\int_{B_r}\vert \vec{\nabla}\otimes \U\vert^2dx\right)^{\frac{1}{2}}   \left(\int_{B_r}\vert \V^{*}\vert^{q} dx \right)^{\frac{1}{q}}    \left(\int_{B_r}\vert  \U\vert^pdx\right)^{\frac{1}{p}}.
\end{equation} 
In this  estimate we still need  to study  the term   $\ds{\left(\int_{B_r}\vert \V^{*}\vert^{q}dx \right)^{\frac{1}{q}}}$. Recall first that the function $\V^{*}$ is defined as $\V^{*}=\V-\V(0)$ where the function $\V$ is given by  the velocity $\U$ in expression  (\ref{def-V}) and since $\U \in \dot{B}^{\frac{3}{p}-\frac{3}{2},\infty}_{\infty}(\Rt)$ then always by expression (\ref{def-V}) we have $\V\in \dot{B}^{\frac{3}{p}-\frac{1}{2},\infty}_{\infty}(\Rt)$. But, recall also that the parameter $p$  verifies $\frac{9}{2}<p<6$ and  then  we have $0<\frac{3}{p}-\frac{1}{2}<\frac{1}{6}$. Thus,  since $\V\in \dot{B}^{\frac{3}{p}-\frac{1}{2},\infty}_{\infty}(\Rt)$ then this function  is an $\alpha-$H\"older continuous function with  $\alpha = \frac{3}{p}-\frac{1}{2}$; and then we can write  $ \ds{\sup_{0<\vert x \vert<r}\frac{\vert \V(x)-\V(0)\vert} {\vert x \vert^{\frac{3}{p}-\frac{1}{2}}}\leq \Vert \V \Vert_{\dot{B}^{\frac{3}{p}-\frac{1}{2},\infty}_{\infty}}}$. \\
\\
With this information and the identity $\V^{*}=\V-\V(0)$, we get back  to the term  $\ds{\left(\int_{B_r}\vert \V^{*}\vert^{q}dx \right)^{\frac{1}{q}}}$ and we write 
$$ \left(\int_{B_r}\vert \V^{*}\vert^{q}dx \right)^{\frac{1}{q}} \leq \left(  \Vert V-\V(0)\Vert_{L^{\infty}(B_r)}\right)    r^{\frac{3}{q}}\leq \left(  r^{\frac{3}{p}-\frac{1}{2}}   \Vert \V \Vert_{\dot{B}^{\frac{3}{p}-\frac{1}{2}, \infty}_{\infty}} \right)  r^{\frac{3}{q}}.$$
But, by the relation  $\frac{1}{2}=\frac{1}{q}+\frac{1}{p}$  we have the identity  $\frac{3}{p}+\frac{3}{q}-\frac{1}{2}=1$, and thus we can write 
$$  \left(\int_{B_r}\vert \V^{*}\vert^{q}dx \right)^{\frac{1}{q}} \leq r  \Vert \V \Vert_{\dot{B}^{\frac{3}{p}-\frac{1}{2}, \infty}_{\infty}}.$$ Moreover,  by equation (\ref{def-V}) we have $\ds{\Vert \V \Vert_{\dot{B}^{\frac{3}{p}-\frac{1}{2}, \infty}_{\infty}} \leq c \Vert \U \Vert_{\dot{B}^{\frac{3}{p}-\frac{3}{2}, \infty}_{\infty}}}$, and since $r<R$ then we  write 
$$ \left(\int_{B_r}\vert \V^{*}\vert^{q}dx \right)^{\frac{1}{q}} \leq c \, R   \Vert \U \Vert_{\dot{B}^{\frac{3}{p}-\frac{3}{2}, \infty}_{\infty}}.$$ 
With this estimate, we get back to inequality (\ref{iden-aux}) and we can write 
$$\ds{\vert  I_3 \vert \leq   c\frac{R}{r-\rho}   \Vert \U \Vert_{\dot{B}^{\frac{3}{p}-\frac{3}{2}, \infty}_{\infty}}    \left(\int_{B_r}\vert \vec{\nabla}\otimes \U\vert^2dx\right)^{\frac{1}{2}} \left(\int_{B_r}\vert  \U\vert^pdx\right)^{\frac{1}{p}}},$$ which is the estimate (\ref{estim-I_3}).  \finpv

\end{document}